\documentclass[10pt,reqno]{article}

\usepackage{amsmath,amsthm}
\usepackage{amsfonts}
\usepackage{amssymb}

\usepackage{graphicx,subfigure}

\usepackage{hyperref}

\usepackage{ifthen} 

\provideboolean{shownotes} 
\setboolean{shownotes}{true} 
\newcommand{\margnote}[1]{
\ifthenelse{\boolean{shownotes}}%
{\marginpar{\raggedright\tiny\texttt{#1}}}%
{}%
}

\newcommand{\hole}[1]{
\ifthenelse{\boolean{shownotes}}%
{\begin{center} \fbox{ \rule {.25cm}{0cm}
\rule[-.1cm]{0cm}{.4cm} \parbox{.85\textwidth}{\begin{center}
\texttt{#1}\end{center}} \rule {.25cm}{0cm}}\end{center}}
{}
}

\theoremstyle{plain}

\theoremstyle{definition}
\newtheorem*{remarks}{Remark}

\theoremstyle{remark}

\newcommand{\R}{\mathbb{R}}

\begin{document}

\title{A chemotactic model for interaction of antagonistic microf\/lora colonies: front asymptotics and numerical simulations}

\author{\textsc{Carlos M\'alaga}\thanks{Departamento de F\'{\i}sica, Facultad
    de Ciencias, Universidad Nacional Aut\'onoma de M\'exico, Circuito
    Exterior s/n, Ciudad Universitaria, C.P. 04510 M\'exico
    D.F. (M\'exico). E-mail: \texttt{cmi@fciencias.unam.mx}} \and \textsc{Antonmar\'{\i}a Minzoni}\thanks{Departamento de Matem\'aticas y Mec\'anica, IIMAS-FENOMEC, Universidad Nacional Aut\'onoma de M\'exico, Apdo. Postal 20-726, C.P. 01000 M\'exico D.F. (Mexico). E-mail: \texttt{tim@mym.iimas.unam.mx; plaza@mym.iimas.unam.mx}} 
  \and \textsc{Ram\'on G. Plaza}$^\dagger$ \and \textsc{Chiara Simeoni}\thanks{Dipartimento di Matematica Pura ed Applicata, Universit\`a degli Studi dell'Aquila, via Vetoio (snc), Coppito I-67010, L'Aquila (Italy). E-mail: \texttt{chiara.simeoni@dm.univaq.it}} }

%
%
%
%
%
%
%
%
%
%
%
%
%

\date{October 25, 2011. Revised: April 24, 2012.}


\maketitle

\begin{abstract}
This paper studies a two-dimensional chemotactic model for two species in which one of them produces a chemo-repellent for the other. It is shown asymptotically and numerically how the chemical inhibits the invasion of a moving front for the second species and how stable steady states, which depend on the chemical concentration, can be reached. The results qualitatively explain experimental observations by Swain and Ray (Microbiol. Res. 164(2), 2009), where colonies of bacteria produce metabolite agents which prevent the invasion of fungi.
\end{abstract}


\section{Introduction}

The movement of biological cells or microscopic organisms in response to chemical gradients is known as \textit{chemotaxis} \cite{Adler1966,MurII3ed}. It pertains to many biological phenomena, such as the motility of bacteria toward certain chemicals \cite{WaAr}, the response of fungal zoospores in the presence of metabolites \cite{IslTah2001}, or the movement of endothelial cells toward angiogenic factors released by tumors \cite{ChSt1}, among others. Since the seminal work by Keller and Segel \cite{KeSe1,KeSe2}, mathematical modeling through systems of partial differential equations has established itself as an efficient tool to describe the macroscopic self-organization patterns of cells occurring in chemotaxis (see \cite{HiPa1} for a review).

This paper studies a chemotactic mathematical model describing the dynamics of two species (or colonies) of cells in which one of them releases a chemo-repellent for the other. Its design and study were motivated by the experimental results of Swain and Ray \cite{SwaRay}, who reported the beneficial activities (known as \textit{biocontrol}) of the bacteria \textit{Bacillus subtilis} in the presence of several harmful microflora found in cowdung, such as the phytopathogenic fungus \textit{Fusarium oxysporum}. In one of the experiments, a uniform concentration of the \textit{F. oxysporum} in a Petri dish was inoculated with the bacteria in two symmetrically located points along the diameter of the dish. The authors observed the inhibition of the in vitro growth of the fungus and the emergence of isolated patterns or strains, free of \textit{F. oxysporum} concentration, which were localized near the places where the \textit{B. subtilis} was applied. These patterns had the form of two circular fronts (see, for example, Figure 2 in \cite{SwaRay}). The bacteria suppressed the growth of the fungus only after 48 hrs. of incubation and the isolation occurred in a stable manner for more than six days (144 hrs.); in other words, the fronts got stabilized and persisted after a long time, suggesting the emergence of stable transient states.
%

It is a well-known fact \cite{LMBFJ,NaBiOb} that the \textit{B. subtilis} produces antifungal metabolites (such as \textit{mycosubtilin}) endowed with biocontrol properties against fungi like the \textit{F. oxysporum} (see, for example, Lecl\'ere \textit{et al.} \cite{LBAGWOTGCJ}). The conditions under which the bacteria and the fungus isolate each other are, however, far from being understood. In order to describe the dynamics of the triggering of fungal suppression by the presence of the metabolite, we propose a chemotactic mechanism of the fungal zoospores moving away from the metabolite gradients. There is strong experimental evidence of zoospore chemotaxis (see, e.g., Islam and Tahara \cite{IslTah2001}) and of \textit{negative} zoospore chemotaxis for certain species (see \cite{AllHar74,CamCar80}). Moreover, it is well-documented that amino acids may act as chemo-repellent molecules for certain zoospores (see Nelson \cite{Nelson06}). Under such considerations (even though it is not clear whether this biocontrol mechanism is of chemotactic nature) we find it reasonable to assume that the inhibition process has a chemotactic component. In addition, up to the authors' knowledge, there are no reliable measurements of diffusion or chemotactic coefficients for the \textit{F. oxysporum} in the literature. 
Thus, the purpose of the model is just to provide a simple mathematical description of the underlying dynamics, which help us to explain qualitatively the emergence of stable invasive fronts for one species and to describe their effects on the biocontrol of the second species through very basic mechanisms. 

As a result, the model studied here consists of a three-component reaction-diffusion system of partial differential equations for the concentrations of the bacteria, the fungus, and the metabolite agent. It is assumed that the latter acts as a chemo-repellent. This fact is expressed through a chemotactic term of Keller-Segel type in the equation for the fungus with the appropriate sign. In addition to diffusion and chemotaxis, we incorporate logistic growth terms for both the bacteria and the fungus, as well as simple linear degradation and production terms for the chemical. Following well-established chemotactic cell-kinetic models (see \cite{HiPa1,WaXue}), we explore the small cell/zoospore diffusivity regime. This means that the diffusion coefficients of both the bacteria and the fungus are small relative to the diffusivity of the chemical. In particular, we found that the solutions reach stationary states in the absence of diffusion for the bacteria. These steady solutions approximate well the metastable patterns observed numerically when the diffusivity coefficient of the bacteria is taken to be very small. Although the metastable states eventually disappear and reach equilibrium solutions, they persist for very long times. (For instance, we conjecture that transient states resembling metastable patterns are what Swain and Ray \cite{SwaRay} observed in their in vitro experiments.) It also to be noted that in both the zero and small bacterial diffusion regimes the chemotactic term prevents the invasion of moving fronts for the fungus concentration, simulating the biocontrol properties of the bacteria. Such fronts get stabilized by the interplay of diffusion, mean curvature and the negative chemotaxis. We show asymptotically and numerically that, for a specific set of parameters, stationary or transient fronts emerge and reach equilibrium configurations. Thus, the simple mathematical model studied here captures the basic dynamical features of the antagonistic interaction of the two species observed in experiments.

It is important to mention that, although the model was designed according to the experiment of Swain and Ray, it is applicable to any system of two antagonistic cell colonies, with one of them producing a chemo-repellent for the other. This paper joins the growing number of works studying multi-species chemotaxis models. Systems of Keller-Segel type with repulsion and attraction of multiple species have been mathematically analyzed by Horstmann in \cite{Hrst1}; the work by Conca \textit{et al.} \cite{CoEsVi} describes the conditions for global existence and blow-up in the case of two-species and one chemo-attractant (see also the study of a one-dimensional attraction-repulsion system by Perthame \textit{et al.} \cite{PSTV}). One-species Keller-Segel systems with logistic growth terms have been studied by Mimura and Tsujikawa \cite{MiTsu1} and by Funaki \textit{et al.} \cite{FuMiT2}. Winkler \cite{Win1,Win3} and Tello and Winkler \cite{TeWi} analyze the circumstances under which logistic growth can limit aggregation effects in one-species models with positive (attractive) chemotaxis. A two-species model (also of Keller-Segel type) with logistic growth has been recently analyzed by Kawaguchi \cite{Kwg}, who studies a three-component system with a chemo-attractant for one of the species. It is to be observed, though, that Kawaguchi reports unstable two-dimensional front-, spot- and wave-like solutions in a chemo-attractive setting, producing branching instabilities. In contrast, this paper provides asymptotic calculations and numerical evidence of the stabilizing effect of chemo-repellent interactions, producing stable/metastable patterns which resemble the ones observed in certain experiments. The system studied here is, we believe, the simplest reasonable chemotactic model capturing the stabilizing effect on invasive fronts of the chemo-repellent interaction of two antagonistic species.

The structure of the paper is as follows. In section \ref{secmodel} the model equations are presented. The derivation of the asymptotic equation of motion for an invasive front is given in section \ref{secasymptotic}. The computation of stationary solutions in the zero-diffusion limit for the bacteria is the content of section \ref{secstationary}. In section \ref{secstability} we discuss the asymptotic conditions under which these stationary solutions are linearly stable. In section \ref{secnumerics} we present the results of our numerical simulations and compare them to the asymptotics of the previous sections. Finally, section \ref{secconclusions} contains a discussion of our results.

\section{Modeling}
\label{secmodel}
We shall denote by $v$ the concentration of the bacteria which produces the chemo-repellent (metabolite) agent. The concentration of the latter will be denoted by $c$. The variable $u$ denotes the concentration of the pathogenic colony (e.g. fungus). In order to describe the antagonistic activities (biocontrol) of these microflora colonies, we propose the following parabolic reaction-diffusion system of equations,
\begin{equation}
\label{modelo}
\begin{aligned}
u_t &= D_u \, \Delta u + \lambda u(u_0-u)(u-u_*) - \nabla \cdot J_c,\\
v_t &= D_v \, \Delta v + \beta v(v_0-v)(v-v_*),\\
c_t &= D_c \,\Delta c + \delta v - \alpha c,  
\end{aligned}
\end{equation}
with $(x,y,t) \in \Omega \times [0,+\infty)$, and where $\Omega \subset \R^2$ is a bounded domain in the plane. According to custom $\Delta = \partial_x^2 + \partial_y^2$ denotes the Laplace operator. The domain $\Omega$ will be taken as a square or a circle with same area, that is
\[
 \Omega = [0,L] \times [0,L] \qquad \mathrm{or} \qquad \Omega = \{x^2 + y^2 \leq L^2/\pi\},
\]  
being $L>0$ the characteristic length of the Petri dish.

The concentration $v$ has diffusivity equal to $D_v > 0$ and the value $v = v_0$ denotes the maximum sustainable population concentration of the bacteria. The production term is of bi-stable (or Nagumo) type \cite{NAY,McK1}, modeling logistic growth with a threshold. This enables the states $v = 0$ and $v = v_0 >0$ to be stable equilibria in the absence of diffusion. We also assume that the unstable equilibrium point $v_*$ satisfies $0 < v_* < v_0$.  Observe that the equation for $v$ is decoupled from the other two. The chemical $c$ is produced by the bacteria at the rate $\delta > 0$ and it degrades at the rate $\alpha > 0$. Its diffusivity coefficient is denoted by $D_c > 0$. Finally, the equation for $u$ is coupled to $v$ and $c$ via a chemotactic term $J_c$. The concentration $u$ diffuses with $D_u > 0$ and the reaction term is also of bi-stable type, with the roots slighted shifted from $u=u_0 >0$ and $u = u_*$ due to the presence of the chemotactic term, as we shall see later on. We assume, of course, that $0 < u_* < u_0$.

The system of equations \eqref{modelo} is further endowed with no-flux boundary conditions of the form
\begin{equation}
\label{bdrycond}
\nabla u \cdot \hat n = 0, \quad \nabla v \cdot \hat n = 0, \quad \nabla c \cdot \hat n = 0, \qquad \mathrm{at} \;\;\;\; \partial \Omega,
\end{equation}
with $\hat n$ being the unit normal to the boundary  $\partial \Omega$, reproducing the physical conditions of the experiments in vitro.


The general form of the chemotactic flux $J_c$ is given by $J_c = - u \chi(c) \nabla c$, where $\chi$ is known as the chemotactic sensitivity function \cite{MurII3ed}. In this work we take
%
\begin{equation}
\label{chemsens}
\chi(c) = \chi_0,
\end{equation}
with $\chi_0 > 0$ constant, which is a customary choice when it is assumed that the species $u$ always responds to a chemosensory stimulus in an uniform manner \cite{MurII3ed,KeSe2}. Moreover, the parameter $\chi_0$ measures the intensity of the chemotaxis. The negative sign in $J_{\mathrm{c}}$ indicates that the metabolite $c$ acts, not as the customary chemo-attractant, but as a chemo-repellent agent blocking the spreading of the species $u$. This type of movement toward lesser concentrations of a chemical is called \textit{negative chemotaxis}. In the case of zoospores, such phenomenon is well-documented \cite{AllHar74}, although it has not been as thoroughly studied as its bacterial counterpart. Therefore, we choose the simplest uniform chemotactic response function \eqref{chemsens}, which is also the first term in the expansion for small concentrations of the more general Lapidus-Schiller kinetic receptor law \cite{LapSch1,ForLauf1}, commonly used for bacteria. Notice that, assuming \eqref{chemsens}, the chemotactic term in the equation for $u$ in \eqref{modelo} has an advection component of the form $\chi_0 \nabla u \cdot \nabla c$ and an apparently less dominant kinetic-type term of the form $\chi_0 u \Delta c$. The effect of the latter is a slight shift on the equilibria when $c$ reaches a stationary state. We observe that, although the system \eqref{ndmodel} models the evolution of three components, it is essentially an scalar model for the fungus concentration $u$, because the equations for $c$ and for $v$ are independent of $u$. In this fashion we capture the simplest mechanism of the inhibition of the fungus due to the presence of the colony $v$ producing a chemo-repellent $c$.

Finally, we build the non-dimensional version of system \eqref{modelo} by making the following substitutions
\[
 D_u \to \frac{D_u}{2D_c}, \quad D_v \to \frac{D_v}{2D_c}, \quad t \to \alpha t, \quad x \to \sqrt{\frac{\alpha}{2D_c}} x, \quad y \to \sqrt{\frac{\alpha}{2D_c}} y,
\]
\[
c \to \frac{c}{c_0}, \quad u \to \frac{u}{u_0}, \quad v \to \frac{v}{v_0}, \quad u_* \to \frac{u_*}{u_0}, \quad v_* \to \frac{v_*}{v_0}, \quad
\]
\[
\lambda \to \lambda \frac{u_0^2}{\alpha}, \quad \beta \to \beta \frac{v_0^2}{\alpha}, \quad \delta \to \delta \frac{v_0}{c_0 \alpha}, \quad \chi_0 \to \chi_0 \frac{c_0}{2D_c}.
\]
The resulting non-dimensional system reads
\begin{equation}
\label{ndmodel}
\begin{aligned}
u_t &= D_u \, \Delta u + \lambda u(1-u)(u-u_*) + \chi_0  \nabla \cdot (u\nabla c),\\
v_t &= D_v \, \Delta v + \beta v(1-v)(v-v_*),\\
c_t &= \tfrac{1}{2} \Delta c + \delta v -  c,  
\end{aligned}
\end{equation}
for $(x,y,t) \in \Omega \times [0,+\infty)$, and where the unstable equilibria satisfy $0 < u_*, v_* <1$. Without loss of generality we assume that $L \sqrt{2D_c/\alpha} = O(1)$, so that the domain will be taken as
\[
 \Omega = [0,1] \times [0,1], \qquad \mathrm{or} \qquad \Omega = \{x^2 + y^2 \leq 1/\pi\},
\]
namely, the unit square or the circle with area equal to one.

In the sequel we analyze solutions to system \eqref{ndmodel}. We are particularly interested in the dynamics of a moving front for the $u$ variable. The asymptotics and stability of such a front, as well as its numerical computation, is the content of the following sections.

\section{Asymptotic equation of motion for the front}
\label{secasymptotic}
When diffusion is small ($D_u \approx 0$), localized patterns or layers in the $u$ variable can be well approximated by interfaces. Here we derive the evolution equation of an interfacial curve when the concentration $u$ is near the equilibrium configuration. Such equation will be controlled by the concentrations $u$ and $c$, and the location of the interface is coupled with $v$ via $c$. We assume that the width of the interface is small. Let us denote
\[
 \Sigma(t) = \{(x,y) \in \Omega \, : \, u(x,y,t) = u_2\}
\]
as the moving front, so that the outer and inner regions are defined by
\[
\begin{array}{l}
\Omega_{\mathrm{in}} = \{ (x,y) \in \Omega \, : \, u(x,y,t) < u_2 \},\\
\Omega_{\mathrm{out}} = \{ (x,y) \in \Omega \, : \, u(x,y,t) > u_2 \}.
\end{array}
\]
The equilibrium point $u_2$ will be determined later. We describe the motion of the interface in local curvilinear coordinates with components $\zeta(x,y,t)$ and $\tau(x,y,t)$, normal and tangent to the interface, respectively, and normalized such that $|\nabla\zeta| = |\nabla \tau| = 1$. When diffusion is small, the dependence of $u$ with respect to the tangent component is negligible and we approximate
\[
 u(x,y,t) \approx \bar u(\zeta(x,y,t)).
\]
Under the assumption that the concentrations $v$ and $c$ have reached stationary (or quasi-stationary) states, we denote their values at the interface as $v_\Sigma$ and $c_\Sigma$, respectively. Upon substitution into the equation for $u$ in \eqref{ndmodel} we obtain an ordinary differential equation for the solution $\bar u$, given by
\begin{equation}
\label{odebaru}
(-s + D_u \kappa - \chi_0 \nabla \zeta \cdot \nabla c_\Sigma)\bar u' = D_u \bar u'' + \lambda \bar u (\bar u - u_*)(1-\bar u) + \chi_0\Delta c_\Sigma \bar u,
\end{equation}
where $s = - \partial_t \zeta$ is the normal velocity, and $\kappa = - \Delta \zeta$ is the local curvature (see, e.g. \cite{OMK}). The dynamics of the front is governed by the time-independent values of $v$ and $c$ at $\Sigma$. 

\subsection{Interface equation of motion}

The total velocity along the normal direction has two main contributions: the velocity of the front in one dimension, which is determined by the nonlinear Nagumo term, and the chemotactic velocity. In order to obtain the former, note that if $s$ and $\kappa$ were independent of $\zeta$ (constant in time) then equation \eqref{odebaru} would resemble a one-dimensional front equation of the form
\begin{equation}
\label{nagumoeqfront}
-s_1 \bar u' = D_u \bar u'' + \lambda \bar u (1-\bar u)(\bar u - u_*) + \chi_0 \Delta c_\Sigma \bar u.
\end{equation}
Since the value of $\Delta c_\Sigma$ does not depend on time (it acts as a coefficient), the reaction term of last equation is also of Nagumo type, resulting in a shift of the roots for the stable/unstable equilibria. Rewriting the right hand side of \eqref{nagumoeqfront}, we obtain
\begin{equation}
\label{nagumoeqshift}
-s_1 \bar u' = D_u \bar u'' + \lambda \bar u (u_1-\bar u)(\bar u - u_2),
\end{equation}
where the equilibrium points $u_1$ and $u_2$ now depend on the value of $\Delta c_\Sigma$ as follows,
\begin{equation}
\label{newspeeds}
\begin{aligned}
u_1 &= \tfrac{1}{2}(1+u_*) + \tfrac{1}{2}\sqrt{(1-u_*)^2 + 4\chi_0 \Delta c_\Sigma/\lambda} \approx 1 + \frac{\chi_0}{\lambda(1-u_*)} {\Delta c_\Sigma},\\
u_2 &= \tfrac{1}{2}(1+u_*) - \tfrac{1}{2}\sqrt{(1-u_*)^2 + 4\chi_0 \Delta c_\Sigma/\lambda} \approx u_* - \frac{\chi_0}{\lambda(1-u_*)} {\Delta c_\Sigma}.
\end{aligned}
\end{equation}

Although the velocity of a one-dimensional front for a general reaction function can be obtained explicitly via a variational characterization \cite{BeDe1}, the Nagumo speed for this modified reaction term is directly computable. Indeed, equation \eqref{nagumoeqshift} is an ordinary differential equation for a monotonic profile $\bar u(\cdot)$ connecting the two stable equilibria $u = 0$ and $u = u_1$. Suppose that $\bar u$ connects these equilibrium points on the left, that is $\bar u (-\infty) = u_1$ and $\bar u(+\infty) = 0$. Therefore the profile is monotonic with $\bar u' < 0$. In the phase space the solution to \eqref{nagumoeqshift} is given by $\bar u' = \phi(\bar u)$, which leads to an equation for $\phi$, namely
\[
 - s_1 \phi = D_u \phi \frac{d\phi}{d \bar u} + \lambda \bar u(u_1-\bar u) (\bar u - u_2),
\]
with boundary conditions $\phi(0) = \phi(u_1) = 0$, and subject to the constraint $\phi < 0$. The solution is of the form $\phi(\bar u) = - \mu \bar u(u_1-\bar u)$ with $\mu > 0$. Since the velocity is independent of $\phi$, upon substitution we obtain $\mu = \sqrt{\lambda / 2D_u}\,$, yielding in turn
\begin{equation}
\label{values1}
s_1 = \sqrt{2\lambda D_u}\Big(\textstyle{\frac{1}{2}} u_1 - u_2\Big).
\end{equation}
The sign of $s_1$ is that of $u_1 - 2u_2 \approx 1 - 2u_* + 3\chi_0 \Delta c_\Sigma  / \lambda (1-u_*)$ in view of \eqref{newspeeds}.

Continued inspection of equation \eqref{odebaru} shows that there is a contribution to the speed of the front due to the chemotaxis. Such term is known as the chemotactic velocity \cite{ForLauf1,KeSe1}. In our two dimensional setting it has the form
\[
s_2 = \chi_0 \nabla \zeta \cdot \nabla c = \chi_0 \frac{dc}{d\zeta}.  
\]

In this fashion, we obtain the interface equation of motion (similar to that of \cite{OMK}, but with a chemotactic contribution):
\begin{equation}
\label{inteqmotion}
s = s_1 - \chi_0 \frac{dc}{d\zeta}  + D_u \kappa.
\end{equation}

\begin{remarks}
We observe that in the absence of chemotaxis the front will stop whenever $s_1 = 0$. The one dimensional case for a single Nagumo equation was studied by Chen \textit{et al.} \cite{ChFlSh} when $s_1$ was taken to be a function of the space variable. It was shown that when $s_1' < 0$ at the steady state $\bar x$ for which $s(\bar x) = 0$ the front is stable. In the present case we do not consider the two dimensional analogue of this situation which arises when the $c$ variables acts as a death factor for the $u$ equation, increasing the death rate of the pathogenic fungus. Instead, we concentrate on the new role of $c$ as a chemorepellent.
\end{remarks}

\subsection{The case of a circular front}

Let us examine the case of a circular front, for which the interfacial curve can be written as
\[
 \Sigma(t) = \{ (x,y) \in \Omega \, : \, u(x,y,t) = u_2 \} = \{ (x,y) \in \Omega \, : \, \sqrt{x^2 + y^2} = R(t) \},
\]
where the function $R : [0,+\infty) \to [0,+\infty)$ indicates the localization and time-evolution of the front. In polar coordinates with $r = \sqrt{x^2 + y^2}$, the normal component is  $\zeta = R(t) - r\,$, so that the curvature is $\kappa = -\Delta \zeta = 1/r$ and the normal velocity is $s = -\partial_t \zeta = -\dot{R}(t)$. Substituting into \eqref{inteqmotion}, the dynamics of a circular front is governed by the equation
\begin{equation}
\label{circlefront}
\dot{R}(t) = - \sqrt{2\lambda D_u} \big({\textstyle{\frac{1}{2}}} u_1 - u_2\big) - \Big(\chi_0 \frac{dc}{dr} + \frac{D_u}{r}\Big)_{|r=R(t)},
\end{equation}
which conveys the combined effects due to curvature, chemotaxis and the Nagumo speed. Here the value of $s_1$ is the one computed in \eqref{values1}, because we are assuming that $\bar u$ connects $u = u_1$ at $r = +\infty$, with $u = 0$ at $r = 0$, in as much as $\zeta$ points out in the direction of the origin and, consequently, $dc/d\zeta = -dc/dr$.

\section{Stationary solutions}
\label{secstationary}
In this section we assume that the concentrations $v$ and $c$ have reached stationary (in the case of zero diffusion $D_v = 0$) or metastable (in the case of small diffusion $0 < D_v \ll 1$) states, which, in addition, have radial symmetry. For simplicity, we suppose that the spatial domain is the circle $\Omega = \{x^2 + y^2 \leq 1/\pi\}$. Since the equation for $v$ is independent from the other two, let us denote by $v_\infty(r)$ the radially-symmetric stationary solution to the equation for $v$ in \eqref{ndmodel}. Therefore, looking for a stationary solution for the chemical concentration $c$ reduces to solving
\begin{equation}
\label{stationaryc}
\begin{aligned}
\tfrac{1}{2} \Delta c + \delta v_\infty -  c &= 0, \;\;\;\; \mathrm{in} \;\;\Omega,\\
\nabla c \cdot \hat n &= 0, \;\;\;\; \mathrm{at} \;\; \partial \Omega.  
\end{aligned}
\end{equation}
We first examine the case of zero diffusion for $v$.

\subsection{The zero diffusion limit}

In the limit of zero diffusion for the bacteria ($D_v = 0$), the equation for $v$  in \eqref{ndmodel} is an ordinary differential equation whose solution tends to the stable equilibrium points $v = 0$ or $v = 1$, depending on the initial spatial distribution. Let us suppose, for example, that the initial condition is a Gaussian of the form
\begin{equation}
\label{gaussianforv}
v(x,y,0) = Ae^{-\omega r^2}, 
\end{equation}
where $r = \sqrt{x^2 + y^2}$ and with $A > 0$, $\omega > 0$. The mass of the initial condition shall be bounded below and above in order to allow the formation of a plateau inside the domain. Hence, we shall further assume that
\begin{equation}
\label{plateaucond}
 1 < \frac{A}{v_*} < e^{\omega/\pi}.
\end{equation}
Therefore, the stationary solution for $v$ is determined by the plateau
\begin{equation}
\label{vinfty}
v_\infty(r) = \begin{cases}
1,& 0 < r < R_1,\\
0,& R_1 < r < 1/\sqrt{\pi},
\end{cases}
\end{equation}
where the radius $R_1$ satisifes $0 < R_1 < 1/\sqrt{\pi}$ and is given by
\begin{equation}
\label{theoR1}
 R_1 = \sqrt{\frac{1}{\omega} \ln (A/v_*)}.
\end{equation}
Then, the solution to \eqref{stationaryc} has the general form
\begin{equation}
\label{solforc}
c(r) =  \begin{cases}
C_1 I_0(\sqrt{2}r) + C_2 K_0(\sqrt{2}r),& R_1 < r < 1/\sqrt{\pi},\\
C_3 I_0(\sqrt{2}r) + \delta, & 0 < r < R_1,
\end{cases}
\end{equation}
where $I_n, K_n$ denote the modified Bessel functions for each $n = 0,1, \ldots$ (see \cite{ArfWeb6e,Bow58}), and $C_i$, $i =1,2,3$, are constants. Recall that $K_0$ diverges at zero. The solution is subject to the Neumann boundary condition at $r = 1/\sqrt{\pi}$ and to $C^1$-matching conditions at $r = R_1$, namely
\[
 \begin{aligned}
  \frac{d}{dr}\Big( C_1 I_0(\sqrt{2}r) + C_2 K_0(\sqrt{2}r) \Big)_{|r=1/\sqrt{\pi}} \; &= 0,\\
  \Big(C_1 I_0(\sqrt{2}r) + C_2 K_0(\sqrt{2}r) - C_3 I_0(\sqrt{2}r) \Big)_{|r=R_1} \; &= \delta,\\
  \frac{d}{dr}\Big( C_1 I_0(\sqrt{2}r) + C_2 K_0(\sqrt{2}r) - C_3 I_0(\sqrt{2}r)\Big)_{|r=R_1} \; &= 0.
 \end{aligned}
\]
%
%
%

Using the known relations (see, e.g., \cite{AbraSteg64}) $I_0'(x) = I_1(x)$, $K_0'(x) = -K_1(x)$ and $K_n(x) I_{n+1}(x) + K_{n+1}(x) I_n(x) = 1/x\,$, we solve for $C_i$ to obtain
\begin{subequations}
\label{lasCs}
\begin{align}
C_1 &= \sqrt{2}\delta R_1 \frac{K_1(\sqrt{2/\pi}) I_1(\sqrt{2}R_1)}{I_1(\sqrt{2/\pi})}, \\
C_2 &= \sqrt{2} \delta R_1 I_1(\sqrt{2}R_1),\\
C_3 &= \frac{\sqrt{2} \delta R_1}{I_1(\sqrt{2/\pi})}\Big(K_1(\sqrt{2/\pi}) I_1(\sqrt{2}R_1) - I_1(\sqrt{2/\pi}) K_1(\sqrt{2}R_1)\Big),
\end{align}
\end{subequations}
completing the form of the solution \eqref{solforc}. 

\subsection{Small diffusion: metastable patterns}

When $D_v > 0$, the solutions to the equation for $v$ in \eqref{ndmodel} with generic initial data and Neumann boundary conditions converge to the stable equilibrium solutions $v = 1$ or $v = 0$ as $t \to +\infty$ (see \cite{CaHo1,CaHo2}). When diffusion is small ($0 < D_v \ll 1$), however, this convergence is very slow and the solutions can exhibit dynamic metastability \cite{CaPe1,CaPe2}. After a short transition time dominated by the Nagumo term, the solution $v$ forms a pattern of layers which is apparently stable. These slowly evolving metastable solutions are neither local minimizers of the energy nor necessarily close to an unstable equilibrium solution. Although these solutions eventually decay to the patternless equilibria when $t \to +\infty$, the time scale for substantial motion of the layers to occur increases exponentially depending on the size of the domain, the diffusion coefficient, and the size of the Nagumo term. For example, in the case of the bistable reaction-diffusion equation in one dimension, namely $u_t = D u_{xx} + f(u)$, Carr and Pego \cite{CaPe1} estimated this mean transient time as $m(T) \sim \exp(Cl/\sqrt{D})$, where $C = \min \{f''(0), f''(1)\}$, and $l > 0$ is of order of the minimum distance between layers and the distance between the layers and the boundary. In our setting, $C = O(\beta)$ and $l = O(1)$, meaning $m(T) \sim O(\exp (\beta/\sqrt{D}))$. 
Even in one dimension, numerical estimations of the transient time seem to depend on $D$, the size of the domain, and the discretization mesh \cite{Est1}. 

Up to the authors' knowledge, there is no analytical estimation available for the mean duration of transient patterns for a bistable reaction-diffusion equation in a two-dimensional spatial domain\footnote{Numerical estimations for a square domain, however, have been provided by Horikawa \cite{Hori1} in the case $D = 1$, showing that $\log_{10} m(T) \sim O(g(l))$, where $g = g(l)$ is a linear function of $l>0$, and $l$ is of order the size of the domain (see figure 6 in \cite{Hori1}).}. Although we cannot provide an estimate for the life span of metastable solutions, our numerical simulations show the existence of a transient state in the $v$ variable that induces, in turn, metastable patterns for the solutions of $c$ and $u$. The stationary solutions for $v$ and $c$ of the previous section in the zero diffusion limit approximate well these metastable patterns (we refer the reader to section \ref{secnumerics} for details).

\subsection{Approximate equilibrium solutions}

We discuss the conditions for the emergence of an equilibrium circular front in the variable $u$ and provide an asymptotic formula to track its location. In this section we suppose that the diffusion coefficient for the species $u$ is small, $D_u \ll 1$, so that the interface equation of motion \eqref{circlefront} is valid in the geometric front propagation limit. As before, the diffussivity coefficient $D_v$ is either zero or very small, $D_v \ll 1$ (no relation between $D_v$ and $D_u$ is assumed at this point). Let us assume that $r = R_0$ is the equilibrium interface position, with $0 < R_0 < 1/\pi$. Here we suppose that either $R(t) \to R_0$ when $t \to +\infty$ (stationary case in the zero-diffusion limit for the bacteria, $D_v = 0$), or that there exist $T_0 = O(1/\lambda) >0$, $T_1 \gg T_0$, and an uniform $\epsilon > 0$, such that
\[
 |R(t) - R_0| < \epsilon, \qquad \mathrm{for} \;\;\; t \in [T_0,T_1],
\]
in the metastable case when diffusion is small $0 < D_v \ll 1$. In both situations, the interface equation for motion \eqref{circlefront} implies the following necessary condition for equilibrium,
\begin{equation}
\label{meracond}
\dot{R}(t)_{|r=R_0} = -\sqrt{2\lambda D_u} \big({\textstyle{\frac{1}{2}}} u_1 - u_2\big) - \frac{D_u}{R_0} - \chi_0 c'(R_0) = \, 0.
\end{equation}

In view that for small times metastable patterns resemble the stationary states in the zero-diffusion limit, we shall use the stationary solutions \eqref{vinfty} and \eqref{solforc}. Therefore, the stationary value of $\Delta c_\Sigma$ at equilibrium is given by $2c(R_0) - \delta v_\infty(R_0)$. After making $\omega$ large enough in the initial distribution \eqref{gaussianforv} for $v$, we assume, without loss of generality, that the equilibrium radius is outside the plateau, i.e. that $R_1 < R_0$. Thus, $v_\infty(R_0) = 0$ and we approximate
\begin{equation}
\label{capprox}
 {\Delta c_\Sigma}_{|r=R_0} \approx 2c(R_0).
\end{equation}
Hence, we approximate the shifted equilibrium points $u_1$ and $u_2$ at the equilibrium position using \eqref{newspeeds}, that yields
\[
 u_1 \approx 1 + \frac{2\chi_0}{\lambda(1-u_*)} c(R_0), \qquad u_2 \approx u_* - \frac{2\chi_0}{\lambda(1-u_*)} c(R_0).
\]

From \eqref{solforc} and \eqref{lasCs}, we have that $c(R_0) = C_1 I_0(\sqrt{2}R_0) + C_2 K_0(\sqrt{2}R_0)$, and $c'(R_0) = \sqrt{2}(C_1 I_1(\sqrt{2}R_0) - C_2 K_1(\sqrt{2}R_0))$. In view that $R_0$ is bounded below, we use the approximations \cite{ArfWeb6e} $I_0(x) \approx 1 + x^2/4$, $I_1(x) \approx x/2$, $K_0(x) \approx \ln(2/x) - \tilde\epsilon$ and $K_1(x) \approx 1/x$ (where $\tilde \epsilon \approx 0.5772$ is the Euler constant), in order to estimate 
\[
\begin{aligned}
c(R_0) &\approx C_1\big(1+{\tfrac{1}{2}}R_0^2 \big) + C_2(\ln\sqrt{2} - \tilde \epsilon - \ln R_0), \\
c'(R_0) &\approx C_1 R_0 - C_2/R_0.
\end{aligned}
\]
Therefore,
\[
 u_2 - \tfrac{1}{2} u_1 \approx u_* - \frac{1}{2}  - \frac{3\chi_0}{\lambda(1-u_*)}\big( C_1\big(1+{\scriptstyle{\frac{1}{2}}}R_0^2 \big) + C_2(\ln\sqrt{2} - \tilde \epsilon - \ln R_0)\big). 
\]
Substituting into \eqref{meracond}, and multiplying by $R_0$, we obtain
\begin{equation}
\label{ptheoR0}
 p(R_0) = 0,
\end{equation}
where the function $p(\cdot)$ is defined by 
\[
 p(x) := a_3 x^3 + a_2 x^2 + a_1 x + a_0 + b x \ln x,
\]
with coefficients
\[
\begin{aligned}
a_3 &= \sqrt{\lambda D_u} \frac{3\chi_0 C_1}{\sqrt{2}\lambda(1-u_*)} = \frac{3\chi_0 C_1\sqrt{D_u}}{\sqrt{2\lambda}(1-u_*)},\\
a_2 &= \chi_0 C_1,\\
a_1 &= \sqrt{2\lambda D_u} \Big( \frac{1}{2} - u_* + \frac{3\chi_0}{\lambda(1-u_*)} \big( C_1 + C_2(\ln \sqrt{2} - \tilde \epsilon\big)  \Big),\\
a_0 &= D_u - \chi_0 C_2,\\
b &= - \sqrt{2 \lambda D_u} \frac{3\chi_0 C_2}{\lambda(1-u_*)} = - \frac{3 \chi_0 C_2 \sqrt{2D_u}}{\sqrt{\lambda}(1-u_*)}.  
 \end{aligned}
\]

Equation \eqref{ptheoR0} provides an asymptotic formula for the equilibrium position. One may apply Newton-Raphson method to compute the zeroes of $p$ in the interval $[0,1/\sqrt{\pi}]$. Once the approximated equilibrium position $R_0$ has been computed numerically, one may use the asymptotic expansions $I_0''(x) \approx \frac{1}{2}(1 + 3x^2/8)$ and $K_0''(x) \approx 1/x^2$, to estimate $c''(R_0) \approx C_1(1+3R_0^2/4) +C_2/R_0^2$. This yields,
\begin{equation}
\label{crucial}
\frac{D_u}{R_0^2} - \chi_0 c''(R_0) \approx \frac{1}{R_0^2} \big( D_u - \chi_0 C_2 \big) - \chi_0 C_1 \big( 1 + \tfrac{3}{4} R_0^2\big).
\end{equation}
The sign of the expression in \eqref{crucial} plays an important role in the stability of the front as we shall see in the next section.

We have thus constructed a steady state solution where the concentration for the variable $u$ equals the equilibrium $u = u_1 \approx 1$ outside a circular region of radius $R_0$ determined by the solutions of the geometric front propagation equation. The chemical $c$ is concentrated inside the plateau of radius $R_1$ and its gradient balances the invading velocities. We will see how this basic state and its stability can explain the behavior of the solutions under biologically relevant initial conditions.

\section{Linearized stability} 
\label{secstability}

To study the stability of the invasive front we consider again equation \eqref{meracond}. Once again we suppose that both $D_u$ and $D_v$ are small, and $D_v$ possibly zero, with no relation between the two. The steady state solution outlined in the previous section can be perturbed in the form $R_0 + \eta(t)$, with $\eta(t) \ll R_0$. Linearization of equation \eqref{meracond} gives
\begin{equation}
\label{starone}
\dot{\eta}(t) = \Big( \frac{D_u}{R_0^2} - \chi_0 c''(R_0) \Big) \eta. 
\end{equation}
This equation shows that the equilibrium point is linearly stable relative to radial perturbations of the front provided that
\begin{equation}
\label{stabilitycond}
\chi_0 c''(R_0) > \frac{D_u}{R_0^2}.
\end{equation}
If $c''(R_0)$ now becomes a slowly evolving function of time, then it follows that, as $c''(R_0)$ diminishes, the equilibrium point will lose stability and the front will start propagating. This will be the case for slowly diffusing chemo-repellent bacteria. As they diffuse with $0 < D_v \ll 1$ the concentration plateau will expand and the equation for $c$ will have a flatter solution as the plateau moves. This will cause a smaller concentration gradient until the final steady state is homogeneous. This prediction is, in fact, verified by the numerical solution which will be discussed in the following section.

The next question is the stability relative to azimuthally dependent perturbations. Again we consider the limit of a moving front and use the equation of motion driven by mean curvature for the interface. If the interface is parametrized by $X(\theta,t) = (x(\theta,t),y(\theta,t))$, where $\theta$ is the polar angle, then the equation of motion for the front takes the form 
\cite{Fif2}
\begin{equation}
\label{startwo}
X_t = \big(c - D_u \kappa(\theta,t) - \nabla c(X(\theta,t)) \cdot \hat n \big) \hat n, 
\end{equation}
where $\hat n$ is the outward unit normal vector, and
\[
\kappa(\theta,t) = \frac{x_\theta y_{\theta\theta} - y_\theta x_{\theta\theta}}{(x_\theta^2 + y_\theta^2)^{3/2}}
\]
is the local curvature. The steady states are given by
\begin{equation}
\label{starthree}
c - D_u \kappa(\theta) - \nabla c(X(\theta)) \cdot \hat n = 0, 
\end{equation}
which provides an ordinary differential equation for the curve $X(\theta)$. In the present situation it is equation \eqref{meracond}. To consider the full linearized stability problem, we take $X(\theta,t) = X_0(\theta) + \xi(\theta,t)$, where $X_0(\theta)$ is the solution of the steady state. The linearized equation for the perturbation $\xi(\theta,t)$ takes the form
\begin{equation}
\label{starfour}
\xi_t = \frac{1}{|X_0'|^3} \Big( 3 \kappa_0 \frac{X_0'\cdot \xi'}{|X_0'|^2} - (X_0')^\perp \cdot \xi'' - X_0'' \cdot (\xi')^\perp \Big) \hat n_0 - \chi_0 (\xi', (D^2c) \hat n_0),
\end{equation}
where $\hat n_0$ is the outer normal to the curve $X_0(\theta)$, the $'\,$ denotes derivative with respect to $\theta$, and $(X_0')^\perp$ and $(\xi')^\perp$ are orthogonal to $X_0'$ and $\xi'$, respectively. The matrix $D^2c$ is the Hessian of the concentration at the equilibrium front. Equation \eqref{starfour} provides the parabolic system for the perturbation front $\xi$. These equations can be further simplified by choosing a representation for $\xi$ in the coordinate system formed by the unit tangent $\hat\tau_0$ and the unit normal $\hat n_0$ to the unperturbed front $X_0$. We thus take
\[
 \xi = p(\theta,t) \hat \tau_0(\theta) + q(\theta,t) \hat n_0(\theta).
\]

It follows immediately from \eqref{starfour} that $p_t = 0$ and hence we may take $p = 0$, in as much as perturbations in the tangential direction amount to a change in parametrization. We are thus left in the present case with an equation for $q$ in the form
\begin{equation}
\label{starfive}
q_t = \frac{D_u}{R_0^2} q_{\theta\theta} + \Big( \frac{D_u}{R_0^2} - \chi_0 c''(R_0) \Big) q, 
\end{equation}
with periodic boundary conditions in $\theta$. Equation \eqref{starfive} shows that azimuthally perturbations decay if the linearized stability condition \eqref{stabilitycond} also holds. It is to be noted that, for an arbitrary azimuthally dependent steady state, the stability equation will be a scalar equation with $\theta$-dependent coefficients. The stability will be determined by the spectrum of the associated differential equation. In particular, \eqref{starfive} shows that fronts with sign changing chemotaxis may be stabilized by diffusion. On the other hand, equation \eqref{starfive} may be helpful to explain instability results (see \cite{Kwg}) for a related system with the opposite chemotaxis sign. 

The issue of nonlinear stability will be addressed numerically in the following section.

\section{Numerical results and nonlinear evolution}
\label{secnumerics}
To study the nonlinear evolution of solutions we solve system \eqref{ndmodel} numerically. For that purpose, we implemented an explicit finite-difference Euler scheme in an appropriately scaled square domain $0 \leq x \leq 1$, $0 \leq y \leq 1$, using a grid spacing with $\Delta x = \Delta y = \frac{1}{N-1} \approx 3.9 \times 10^{-3}$ and $N = 256$. The time step was set as $\Delta t = 10^{-3} \times \Delta x \approx 3.9 \times 10^{-6}$. Hence, the Courant number $\mu = \Delta t / (\Delta x)^2$ is approximately $\mu \approx \frac{1}{4}$, assuring numerical stability. In addition, we imposed zero-flux boundary conditions for all species.

Although explicit Euler schemes are not commonly used to solve stiff equations (because they require very small time steps in order to avoid instabilities), they are particularly easily implemented to run in parallel on a Graphics Processing Unit (GPU). In this way, the drawback of a small time step is compensated by the high performance parallel computation capabilities of a GPU with hundreds of processors. Our computations were performed on a commodity-type NVIDIA$^\copyright$ GeForce GTX 480 graphics card with 480 stream processors. In this fashion, we were able to compute tenths of millions of time steps in a few hours.

Finally, to perform the numerical calculations it is necessary to determine the non-dimensional parameter values involved in model \eqref{ndmodel}. Table \ref{tableparam} shows the set of parameter values used during our simulations.

\begin{table}[t]
\caption{Non-dimensional parameter values used in the computations.}
\label{tableparam}
\begin{tabular}{|l|c|r|}\hline \rule[-3mm]{0mm}{8mm}
Description & Symbol & Value \\ \hline \hline
\rule[-1mm]{0mm}{5mm} Diffusion coefficient of $u$ 
 & $D_u$ & 0.01\\
\rule[-1mm]{0mm}{5mm} Diffusion coefficient of $v$ 
 & $D_v$ & $10^{-5}$ or 0\\
\rule[-1mm]{0mm}{5mm} Rate of production of $c$ & $\delta$ & 10.0 \\
\rule[-1mm]{0mm}{5mm} Reaction coefficient of $v$
& $\beta$ & 8.0 \\
\rule[-1mm]{0mm}{5mm} Reaction coefficient of $u$ 
& $\lambda$ & 60.0 \\
\rule[-1mm]{0mm}{5mm} Chemotactic sensitivity & $\chi_0$ & 3.2  \\
\rule[-1mm]{0mm}{5mm} Unstable equilibrium for $u$ & $u_*$ & 0.2 \\
\rule[-3mm]{0mm}{7mm} Unstable equilibrium for $v$ & $v_*$ & 0.5 \\
\hline 
\end{tabular}
\end{table}

\subsection{Zero-diffusion}

We begin by studying system \eqref{ndmodel} with $D_v = 0$. The initial conditions for the bacteria (species $v$) is a colony concentrated at the origin of the domain, corresponding to a centered Gaussian of the form \eqref{gaussianforv}, with $A=3$, $\omega = 1000$ and $v_* = 0.5$, that is
\begin{equation}
\label{civT1}
v(x,y,0) = 3e^{-1000((x-0.5)^2+(y-0.5)^2)}.
\end{equation}
It produces a chemical $c$ which initiates with concentration $c(x,y,0) = 0$. The initial value for the species $u$ is
\begin{equation}
\label{ciuT1}
u(x,y,0) = 10 \Big(e^{1000((x-0.2)^2+(y-0.2)^2)} + e^{-1000((x-0.2)^2+(y-0.2)^2)}\Big)^{-1},
\end{equation}
which represents a localized colony at one of the corners of the domain.

\begin{figure}[ht!]
\begin{center}
        \subfigure[$t = 0.1961$.]{
            \label{figT1v-condini}
            \includegraphics[scale=0.22]{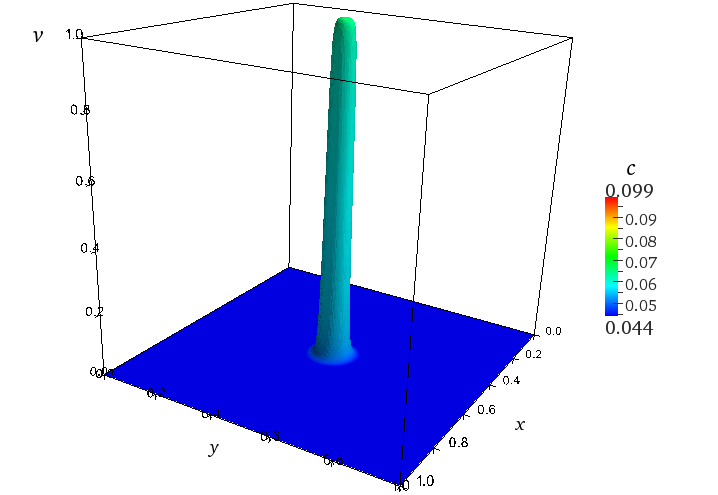}
        }
	\subfigure[$t = 0.9804$.]{
            \label{figT1v-t09804}
            \includegraphics[scale=0.22]{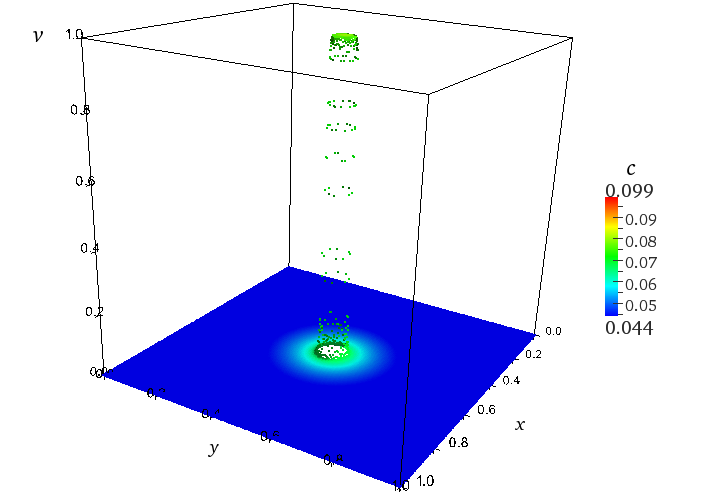}
        }\\
	\subfigure[$t = 1.9608$.]{
            \label{figT1v-t19608}
            \includegraphics[scale=0.22]{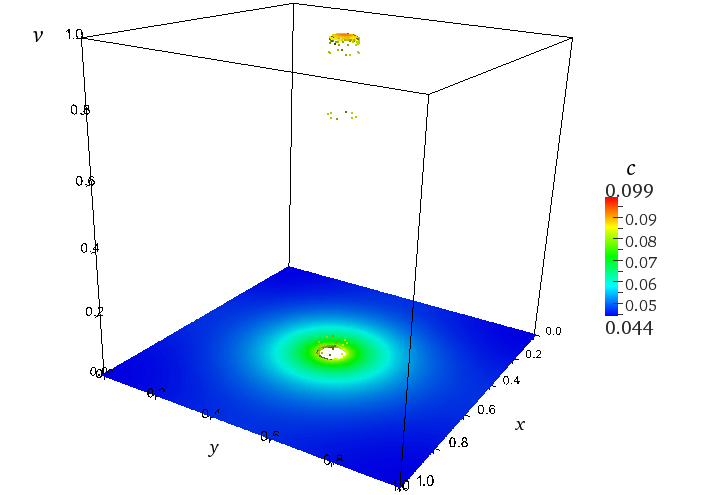}
        }
	\subfigure[$t = 9.8039$.]{
            \label{figT1v-t98039}
            \includegraphics[scale=0.22]{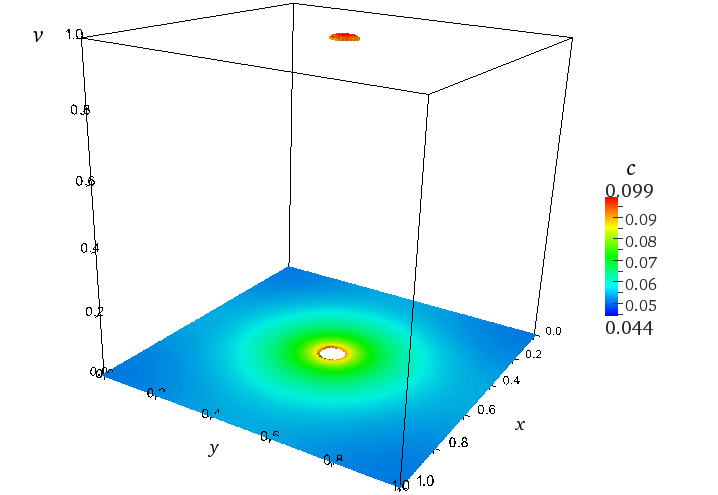}
        }\\
	\subfigure[Chemical concentration $c$ at time $t=9.8039$.]{
	    \label{figT1c-t98039}
	    \includegraphics[scale=0.22]{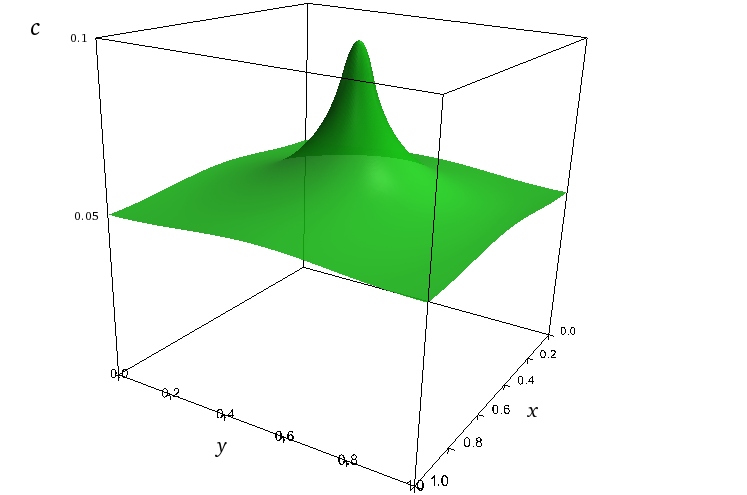}
	}
\end{center}
\caption{\label{figT1v}Figures (a) - (d) show the concentration $v$ of the bacteria colony as time evolves. The bacteria $v$ reaches a stationary state given by the plateau in figure (d). In this figures, the concentration for $v$ is colored according to the concentration of the chemical it produces (see color chart on the right). The steady state for the concentration of the chemical $c$ at time $t = 9.8039$ is shown in figure (e).}
\end{figure}

\begin{figure}[ht!]
\begin{center}
        \subfigure[$t = 0.1961$.]{
            \label{figT1u-condini}
            \includegraphics[scale=0.19]{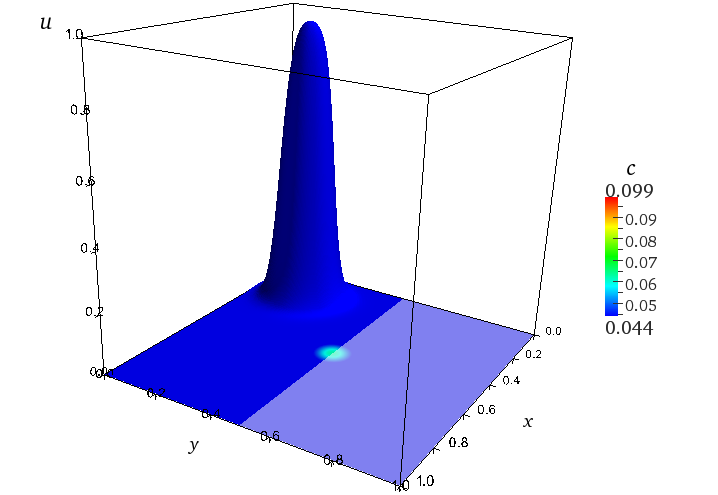}
        }
	\subfigure[$t = 0.7843$.]{
            \label{figT1u-t07843}
            \includegraphics[scale=0.19]{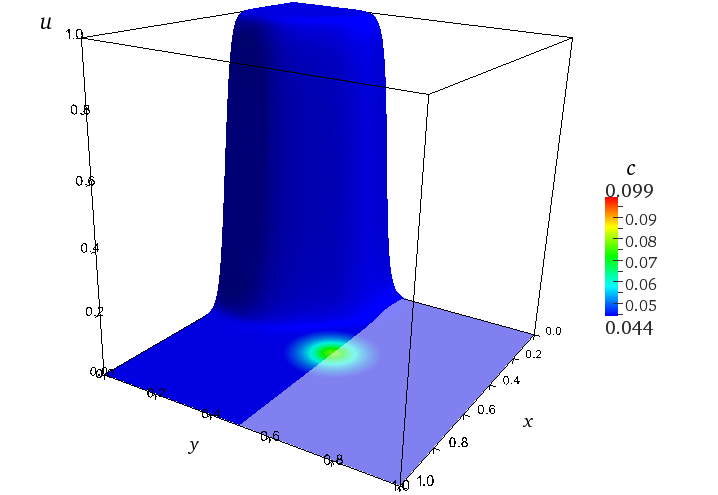}
        }\\
	\subfigure[$t = 1.1765$.]{
            \label{figT1u-t11765}
            \includegraphics[scale=0.19]{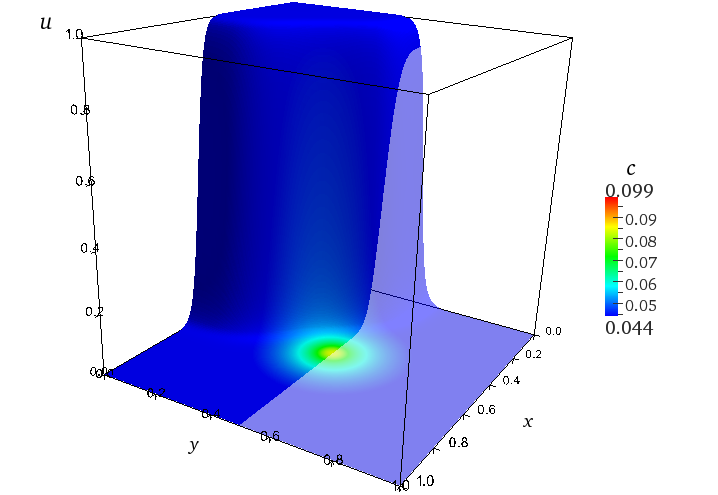}
        }
	\subfigure[$t = 1.7647$.]{
            \label{figT1u-t17647}
            \includegraphics[scale=0.19]{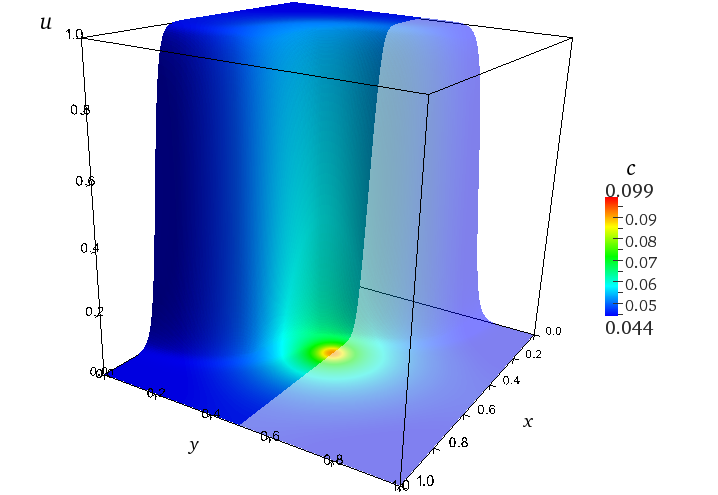}
        }\\
	\subfigure[$t = 2.5490$.]{
            \label{figT1u-t25490}
            \includegraphics[scale=0.19]{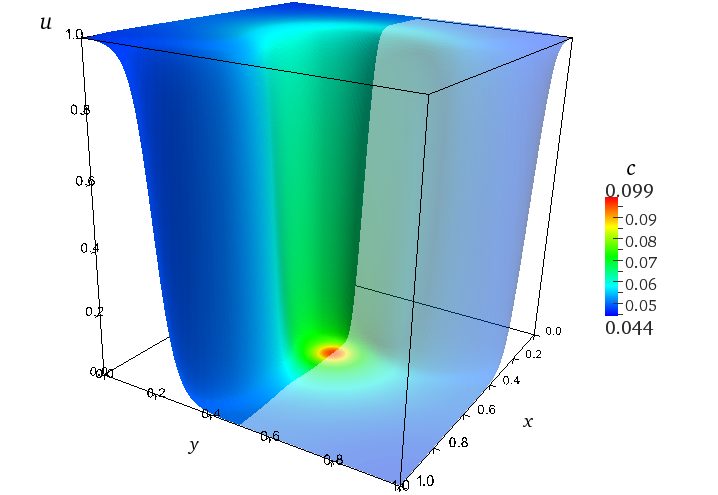}
        }
	\subfigure[$t = 3.1373$.]{
            \label{figT1u-t31373}
            \includegraphics[scale=0.19]{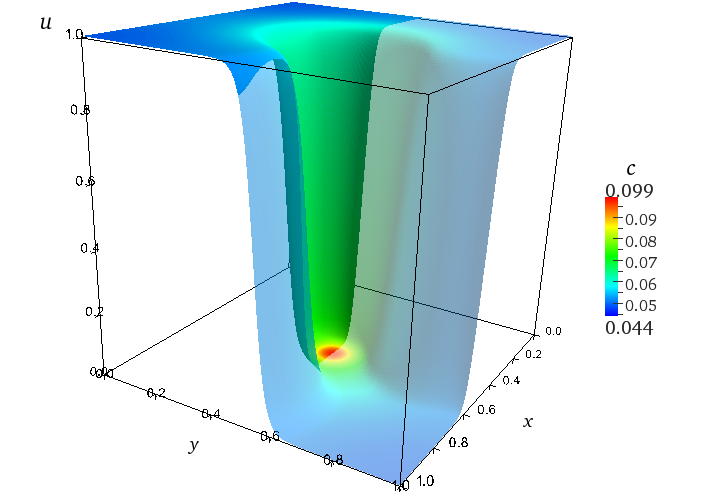}
        }\\
	\subfigure[$t = 3.9216$.]{
            \label{figT1u-t39216}
            \includegraphics[scale=0.19]{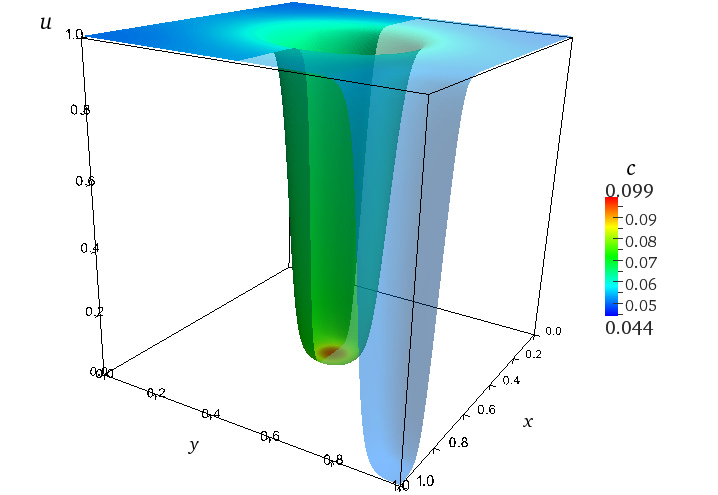}
        }
	\subfigure[$t = 9.8039$.]{
            \label{figT1u-t983039}
            \includegraphics[scale=0.19]{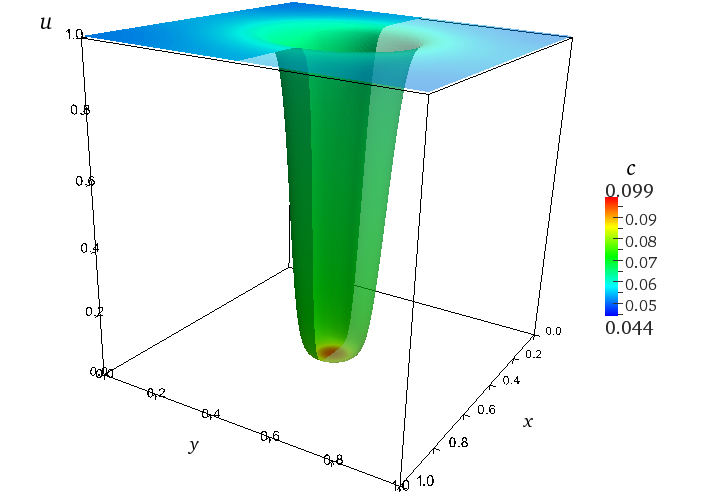}
        }
	\\
\end{center}
\caption{\label{figT1u}Concentration $u$ of the harmful colony (fungus) for different times in the zero-diffusion limit for $v$ (i.e. with $D_v = 0$). Figure (a) shows the concentration of $u$ at time $t =0.1961$, near the initial condition \eqref{ciuT1}. The latter was localized near one of the corners of the domain. Observe that $u$ diffuses, senses the chemo-repellent localized in the center of the domain, and reaches a stationary state depicted in figure (h).
}
\end{figure}

\begin{figure}[ht!]
\begin{center}
\includegraphics[scale=0.45]{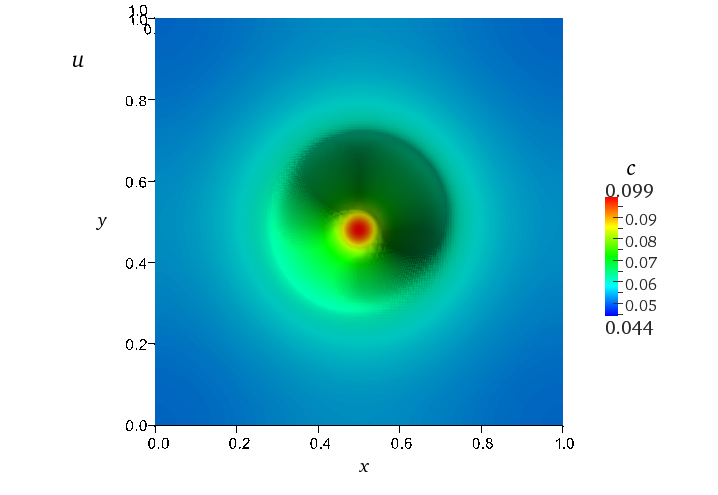}
\end{center}
\caption{\label{figT1utop}Stationary end state for the concentration $u$ of the harmful colony (fungus), at time $t = 9.8039$. It is coloured based on the concentration $c$.}
\end{figure}

The evolution of the solution to system \eqref{ndmodel} is displayed in figures \ref{figT1v} to \ref{figT1utop}. In figures \ref{figT1v-condini}, \ref{figT1v-t09804}, \ref{figT1v-t19608} and \ref{figT1v-t98039} we observe the formation of a sharp plateau for the bacteria concentration $v$. The cross section of the plateau has radius given analytically by \eqref{theoR1} and, with the parameter values $A = 3$, $\omega = 1000$ and $v_* = 0.5$, then it has an approximated value of $R_1 = 0.0423$. (We can also evaluate the constants in \eqref{lasCs} to obtain a theoretical stationary solution $c$ given in \eqref{solforc}). In these numerical simulations, the computed value for the radius of the plateau $R_1$ was approximately $R_1 \approx 0.0412\,$, which gives a relative error of 2.6\% with respect to its theoretical value. The chemical concentration $c$ reaches a steady state which is displayed in figure \ref{figT1c-t98039} at time $t = 9.8039$.

Figures \ref{figT1u-condini} and \ref{figT1u-t07843} display the formation of the initial invasive front in the variable $u$ which generates from the initial condition. Note that it becomes perpendicular to the boundary due to the Neumann conditions. At the time steps of figures \ref{figT1u-t11765} and \ref{figT1u-t17647}, the plateau of bacteria and the steady state for the chemical are already formed, causing the emergence of a repelling region in the center of the domain for the front in the variable $u$. In figures \ref{figT1u-t25490} and \ref{figT1u-t31373} it is shown how the front responds to the chemical gradient by bending according to equation \eqref{startwo}. Finally, the front encircles the repelling region forming the steady state whose depression is centered around the bacteria. The steady state for $u$ is depicted in figure \ref{figT1u-t983039} at time $t = 9.8039$. All surfaces in figures \ref{figT1u-condini} - \ref{figT1u-t983039} are colored according to the chemical concentration for an easier interpretation of the chemotactic effects. 

Finally, figure \ref{figT1utop} shows a top view of the final steady state for the front. Recall that we may compute a theoretical radius $R_0$ of equilibrium by solving equation \eqref{ptheoR0}. With the parameter values depicted in Table \ref{tableparam} and taking $A = 3$, $\omega = 1000$, and $\delta = 10$ in \eqref{gaussianforv}, there is only one zero of $p(\cdot)$ in $[0,1/\pi]$, and hence, the predicted equilibrium radius is approximately
\begin{equation}
\label{theoreticalR0}
R_0 \approx 0.1315.
\end{equation}

It is to be noticed that the computed numerical value for the equilibrium radius of the solution depicted in figure \ref{figT1utop} is $R_0 \approx 0.1235$. A comparison with the theoretical value in \eqref{theoreticalR0} gives a relative error of 6\%. Moreover, in both cases the stability condition \eqref{stabilitycond} holds. For example, for the theoretical value of $R_0$ we obtain
\begin{equation}
\label{compara}
\chi_0 c''(R_0) \approx 3.4409, \qquad \frac{D_u}{R_0^2} \approx 0.5783,
\end{equation}
yielding the negative sign of \eqref{crucial} in order to fulfill the stability condition. It is important to remark that condition \eqref{stabilitycond} is necessary for stability of the front under both radial and azimuthal perturbations.

\subsection{Small diffusion: $0 < D_v \ll 1$}

\begin{figure}[ht!]
\begin{center}
        \subfigure[$t = 0.980$.]{
            \label{v-t0-980}
            \includegraphics[scale=0.22]{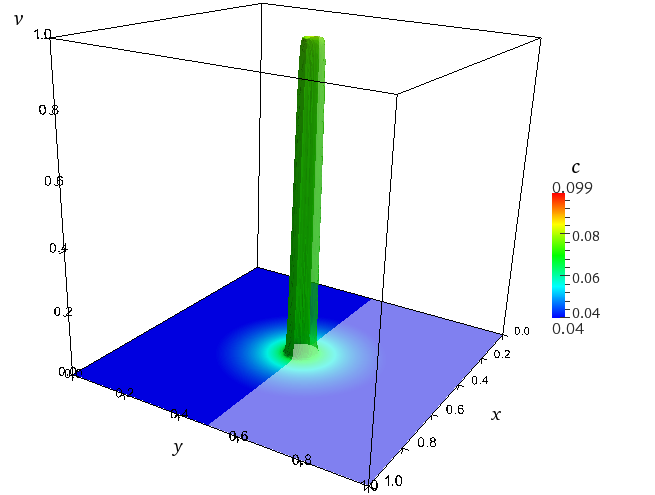}
        }
	\subfigure[$t = 9.804$.]{
            \label{v-t9-804}
            \includegraphics[scale=0.22]{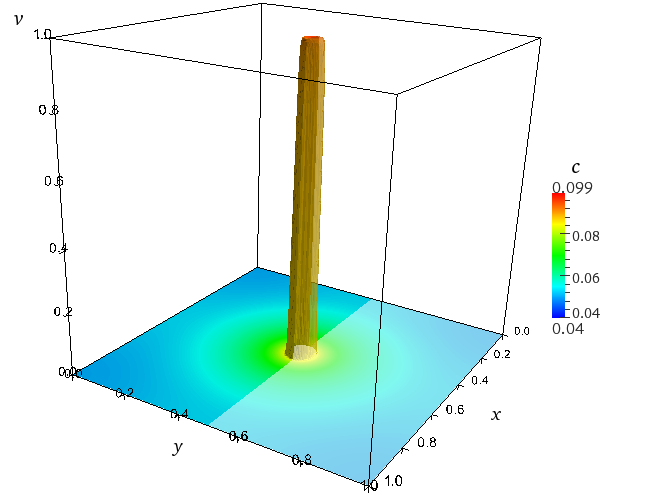}
        }\\
	\subfigure[$t = 99.020$.]{
            \label{v-t99-020}
            \includegraphics[scale=0.22]{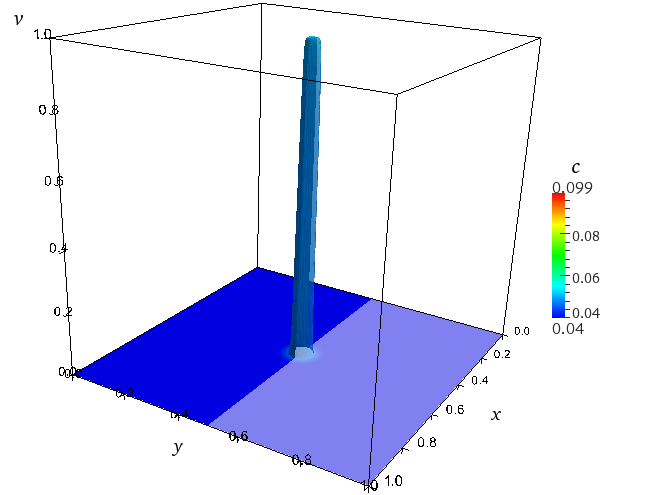}
        }
	\subfigure[$t = 166.667$.]{
            \label{v-t166-667}
            \includegraphics[scale=0.22]{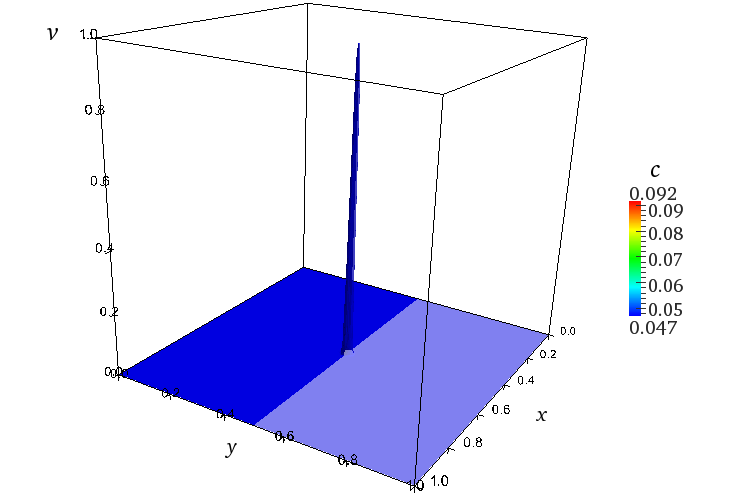}
        }
\end{center}
\caption{\label{figT2v}Figures (a) - (d) show the concentration $v$ of the bacteria colony as time evolves when $D_v = 10^{-5}$. The concentration of $v$ is now a smooth function (due to diffusion) that resembles the plateau of figure \ref{figT1v-t98039} when $D_v = 0$ for short times. Figure (a) depicts this metastable state. The transient, however, tends to the equilibrium point $v = 0$ when $t \to + \infty$. In this figures, the concentration for $v$ is colored according to the concentration of the chemical it produces (see color chart on the right).}
\end{figure}

In order to explore numerically the effects of $D_v \neq 0$, we took $D_v = 10^{-5}$ as the diffusion coefficient for the bacteria, together with the same initial conditions \eqref{civT1} and \eqref{ciuT1} for the bacteria and the fungus, respectively. The initial distribution for the chemical was taken, once again, as $c(x,y,0) = 0$. In such a case, the bacteria will themselves evolve according to a spreading front in the bistable reaction-diffusion equation. It is known \cite{CaHo1,CaHo2} that under Neumann boundary conditions and positive diffusion, the solutions will eventually diffuse and decay to one of the two stable constant equilibrium solutions $v=0$ or $v=1$. Under the parameter values of our simulations and the initial condition considered, we observed the case in which the bacteria concentration eventually vanishes for very large times as it is shown in figure \ref{figT2v}.

After a short time, however, a metastable solution is formed. This pattern is smooth due to diffusion but resembles the numerical plateau of figure \ref{figT1v}. Observe that this metastable state is formed at a time before $t = 0.980$ (figure \ref{v-t0-980}), and persists for a very long time of order $t = 100$ (figure \ref{v-t99-020}). The concentration of bacteria eventually vanishes (figure \ref{v-t166-667}) and, before a time of order $t = 180$, it has already reached the steady state $v = 0$.

\begin{figure}[ht!]
\begin{center}
        \subfigure[$t = 9.8039$.]{
            \label{figT2c_t_9_8039}
            \includegraphics[scale=0.22]{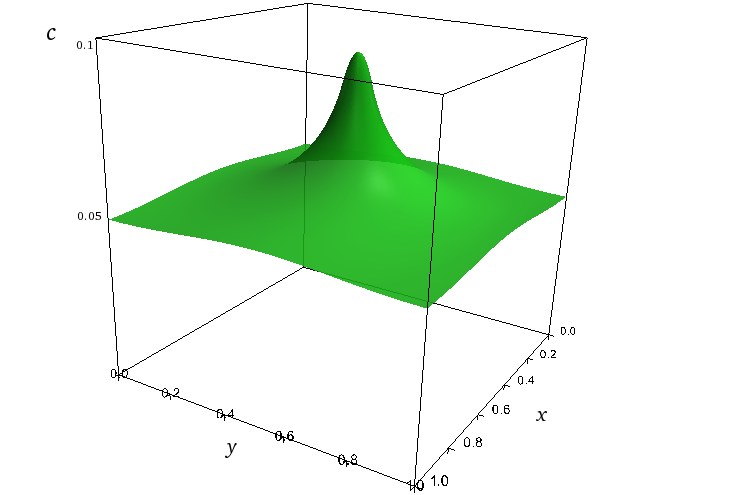}
        }
	\subfigure[$t = 49.020$.]{
            \label{figT2c_t_49_020}
            \includegraphics[scale=0.22]{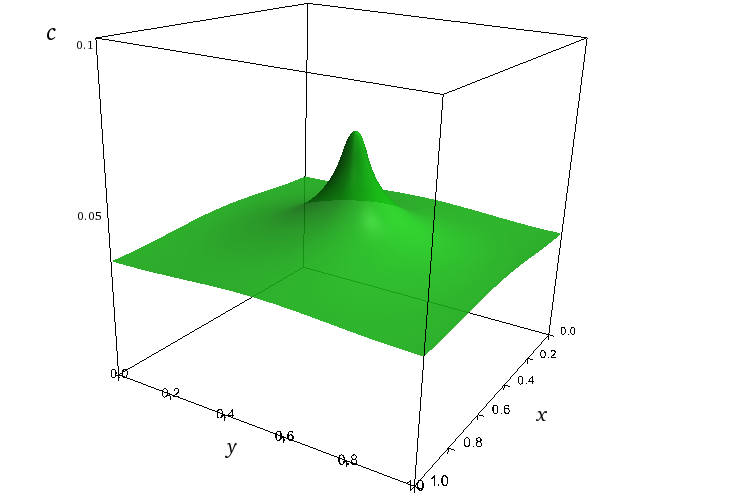}
        }\\
	\subfigure[$t = 88.235$.]{
            \label{figT2c_t_88_235}
            \includegraphics[scale=0.22]{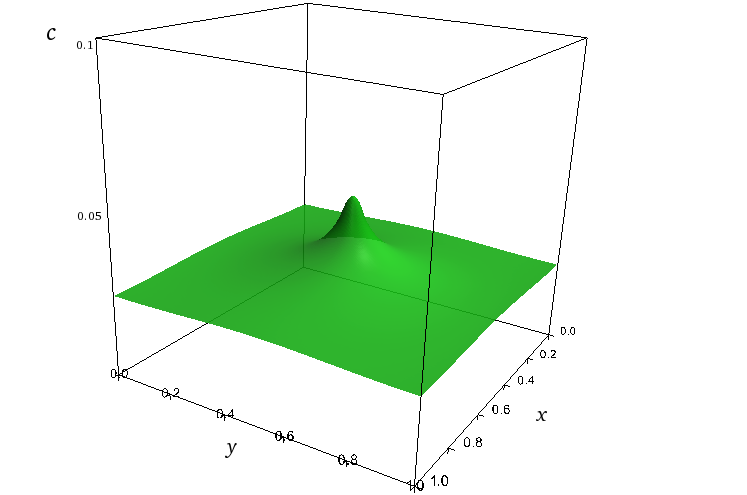}
        }
	\subfigure[$t = 147.159$.]{
            \label{figT2c_t_147_159}
            \includegraphics[scale=0.22]{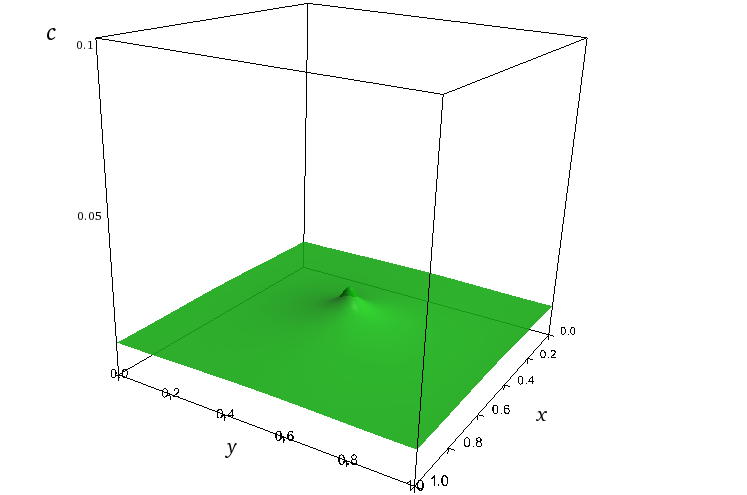}
        }
\end{center}
\caption{\label{figT2c}Figures (a) - (d) show the concentration $c$ of the chemical for different times, when $D_v = 10^{-5}$. Although eventually $c$ reaches the constant steady state $c = 0$, figures (a) and (b) show the metastable state for short times. Note that figure (a) resembles the stationary state computed when $D_v = 0$ in figure \ref{figT1c-t98039}.}
\end{figure}

\begin{figure}[h]
\begin{center}
        \subfigure[$t = 0$.]{
            \label{u-t0}
            \includegraphics[scale=0.19]{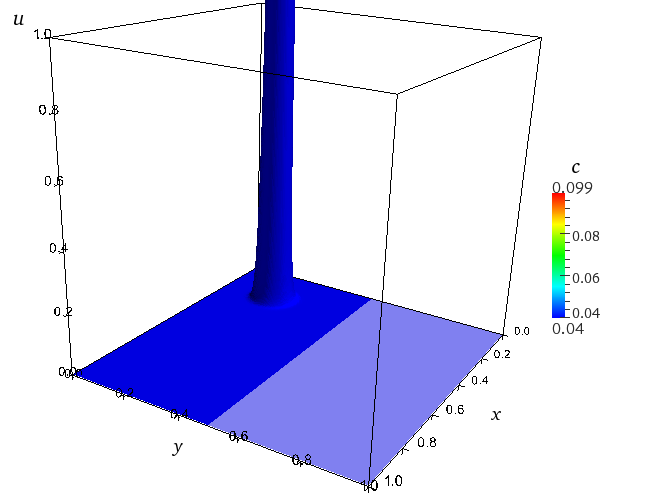}
        }
	\subfigure[$t = 1.961$.]{
            \label{u-t1-961}
            \includegraphics[scale=0.19]{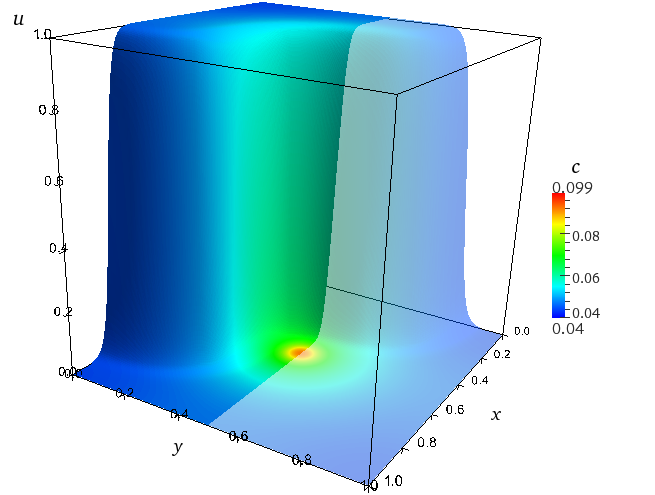}
        }\\
	\subfigure[$t = 4.902$.]{
            \label{u-t4-902}
            \includegraphics[scale=0.19]{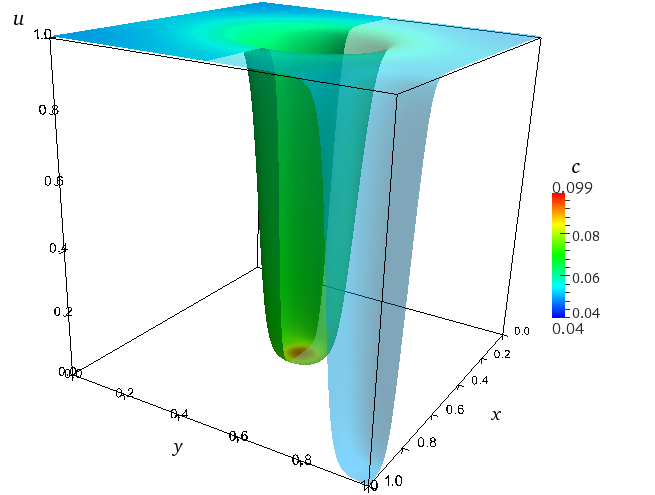}
        }
	\subfigure[$t = 5.882$.]{
            \label{u-t5-882}
            \includegraphics[scale=0.19]{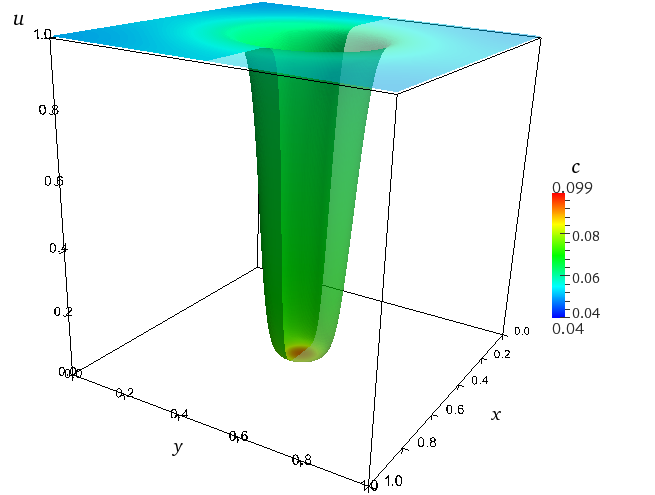}
        }\\
	\subfigure[$t = 50$.]{
            \label{u-t50-000}
            \includegraphics[scale=0.19]{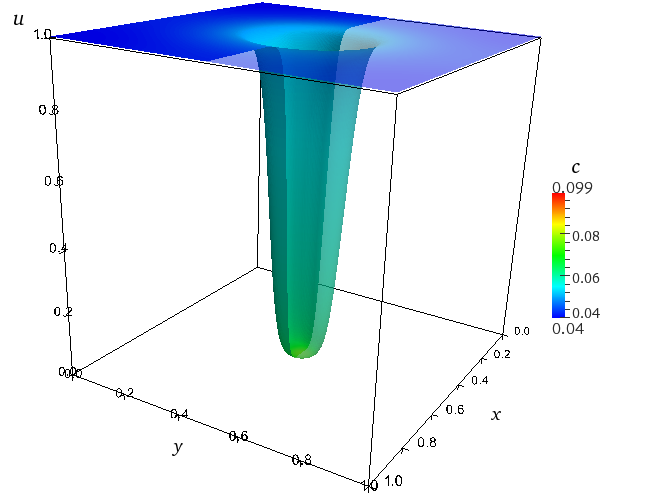}
        }
	\subfigure[$t = 99.020$.]{
            \label{u-t99-020}
            \includegraphics[scale=0.19]{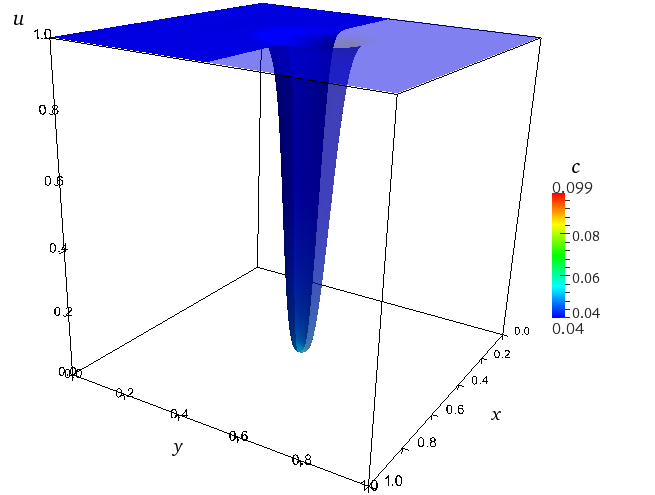}
        }\\
	\subfigure[$t = 147.059$.]{
            \label{u-t147-059}
            \includegraphics[scale=0.19]{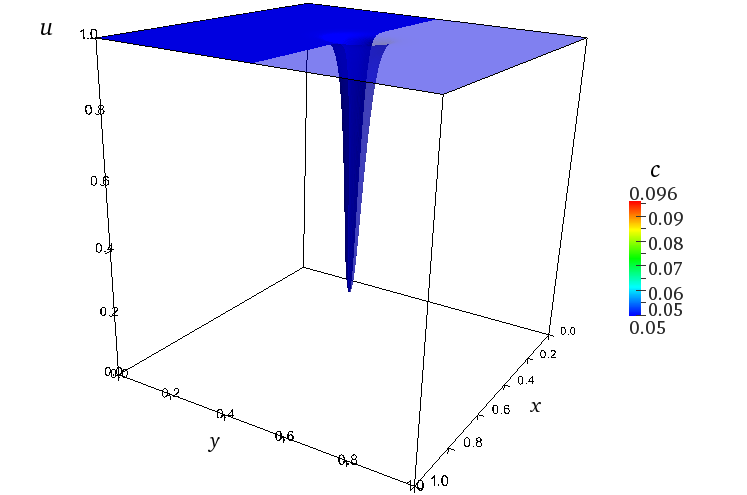}
        }
	\subfigure[$t = 166.667$.]{
            \label{u-t166-667}
            \includegraphics[scale=0.19]{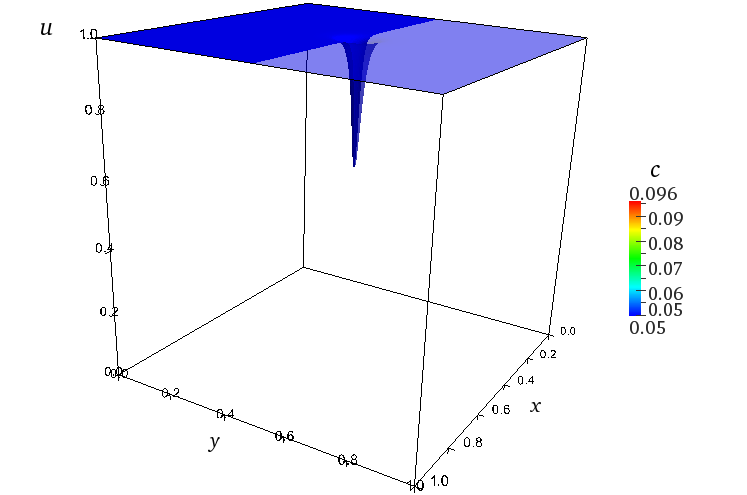}
        }
	\\
\end{center}
\caption{\label{figT2u}Fungus concentration $u$ for different times when $D_v = 10^{-5}$. Observe the emergence of a metastable pattern which resembles the steady state computed with $D_v = 0$ (figure \ref{figT1u-t983039}). This metastable solution changes very little between times $t = 5.882$ (figure (d)) and $t = 50$ (figure (e)). Due to bacterial diffusion, the chemical concentration decays allowing the front to invade all the domain as shown in figures (f), (g) and (h).}
\end{figure}

As this behavior occurs, the concentration $c$ of the chemical also forms a slowly evolving metastable pattern. Figure \ref{figT2c_t_9_8039} depicts this solution, which persists for a very long time. Notice the resemblance of the steady concentration $c$ when $D_v = 0$ of figure \ref{figT1c-t98039} with the concentration when $D_v = 10^{-5}$ for small times shown in figure \ref{figT2c_t_9_8039}. After a long time of ordet $t = 150$, the chemical $c$ reaches the steady state $c = 0$ as it is shown in figures \ref{figT2c_t_49_020}, \ref{figT2c_t_88_235} and \ref{figT2c_t_147_159}. In this case the chemotactic gradient vanishes and, as expected from equation \eqref{starone}, it is shown in figures \ref{u-t0} to \ref{u-t99-020} that the front of the fungus eventually covers the domain since the steady state is unstable. As $D_v$ becomes smaller, the time of existence of the steady (metastable) state is longer.

\begin{figure}[ht!]
\begin{center}
\includegraphics[scale=0.4]{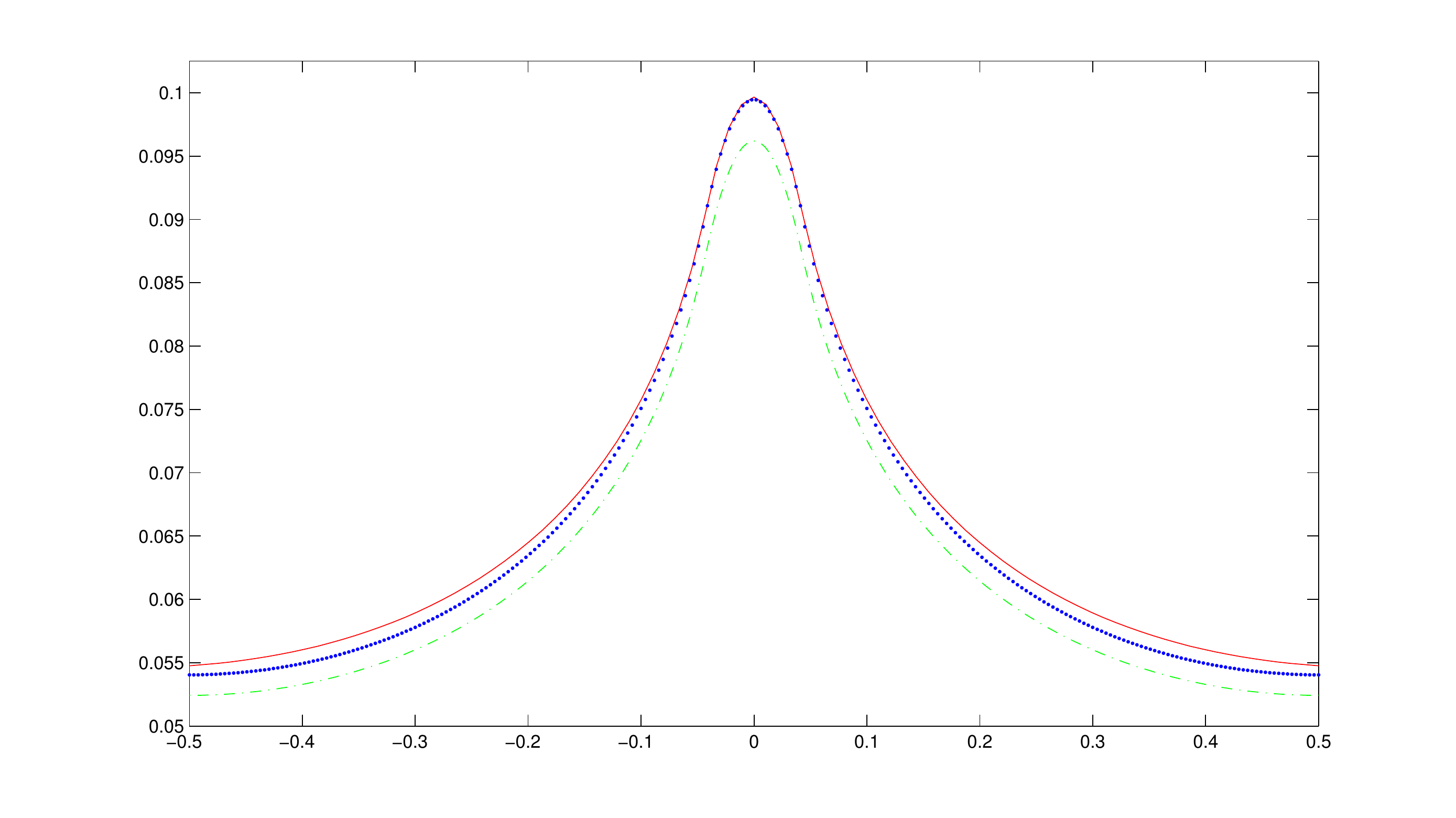}
\end{center}
\caption{\label{figslicec}The continuous (red) plot shows the cross-section for $c$ computed analytically in \eqref{solforc} in the zero diffusion limit $D_v = 0$. The dotted plot (blue) shows the numerical solution for $c$ when $D_v = 0$ at time $t = 9.8039$, when it has reached the stationary state. The dashed line (green) shows the numerically computed concentration $c$ when $D_v = 10^{-5}$, also at time $t = 9.8039$, depicting the metastable state.}
\end{figure}

\begin{figure}[ht!]
\begin{center}
\includegraphics[scale=0.4]{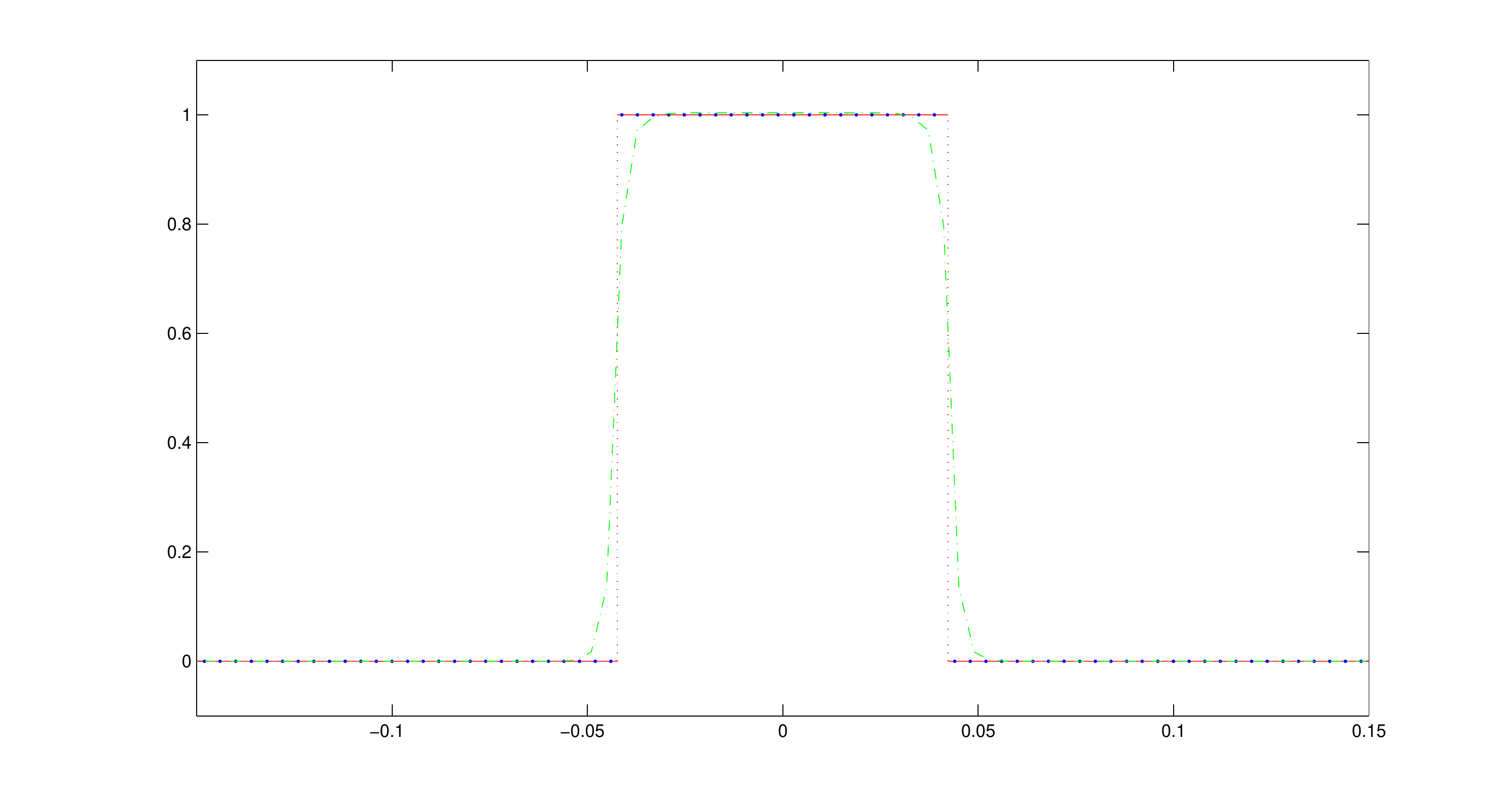}
\end{center}
\caption{\label{figslicev}The continuous (red) plot shows the cross-section of the stationary state (plateau) $v$ computed analytically in \eqref{vinfty} in the zero diffusion limit $D_v = 0$. The dotted plot (blue) shows the numerically computed plateau for $v$ when $D_v = 0$ at time $t = 9.8039$. The dashed line (green) represents the numerically computed concentration $v$ when $D_v = 10^{-5}$, also at time $t = 9.8039$ (metastable state).}
\end{figure}

\begin{figure}[ht!]
\begin{center}
\includegraphics[scale=0.4]{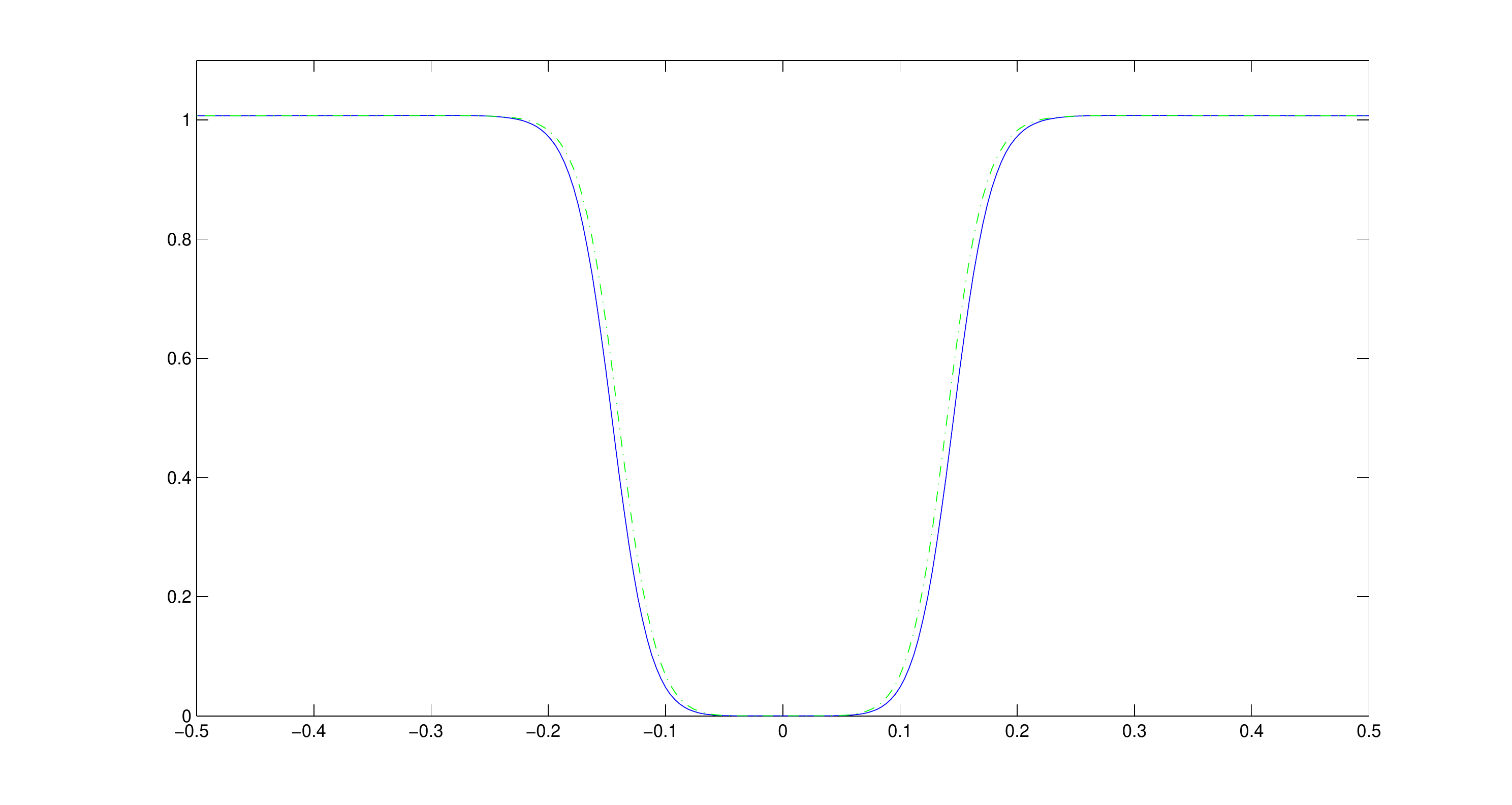}
\end{center}
\caption{\label{figsliceu}The continuous (blue) plot shows the cross-section of the numerically computed stationary profile for the concentration $u$ in the zero diffusion limit $D_v = 0$ at time $t = 9.8039$. The dotted plot (green) represents the numerically computed concentration $u$ in a metastable state at time $t = 9.8039$ when the diffusion coefficient $ 0 < D_v = 10^{-5}$ is small. Observe that the stationary value computed in the zero diffusion limit is a good approximation of the metastable concentration for short times.}
\end{figure}

It is important to remark that the computed solutions in the $D_v = 0$ case are good approximations of the metastable states when we perturb the diffusion parameter $D_v$. Figures \ref{figslicec}, \ref{figslicev} and \ref{figsliceu} show the corresponding cross sections for the final steady states in the zero diffusion limit and the metastable states with $D_v = 10^{-5}$ for short times. For instance, figure \ref{figslicec} shows cross sections of the concentration of the chemical $c$. The continuous (red) plot corresponds to the analytic steady solution \eqref{solforc} computed for a circular domain in the zero diffusion limit $D_v = 0$. The dotted graph (blue) shows the numerically computed solution $c$ on a square domain with $D_v = 0$ at time $t = 9.8039$, when the steady state has been already reached. Notice that it approximates well the analytic stationary solution. The error, which can be observed in the tails of both graphs, is due to the Neumann conditions computed on a circle (analytical solution) and on a square (numerical solution). The dashed (green) graph depicts the concentration $c$ computed numerically with $D_v = 10^{-5}$ at time $t = 9.8039$, when a metastable state has been reached, but in a short time step relatively to the relaxation time with small diffusion. It can be observed that the numerical solutions, both in the zero diffusion and small diffusion limits, approximate well the analytical solution \eqref{solforc}.

In the same fashion, figure \ref{figslicev} shows the cross sections for the concentration $v$ of bacteria. The solid plot (red) shows the plateau \eqref{vinfty} with an analytical radius equal to $R_1 = 0.0423$. The dotted (blue) graph represents the numerically computed plateau with $D_v = 0$ at time $t = 9.8039$, with radius $R_1 = 0.0412$. The dashed graph (green) shows the concentration for $v$ at time $t = 9.8039$ when computed numerically with diffusion equal to $D_v = 10^{-5}$. It is smooth due to diffusion and it corresponds to a metastable state which remains for a long time, until it eventually reaches a constant state $v = 0$.

Figure \ref{figsliceu} shows the cross section for the numerically computed concentrations for $u$ in the zero diffusion case $D_v = 0$ (continuous plot in blue), and in the case of small diffusion $D_v = 10^{-5}$ (dashed plot in green). The former corresponds to a steady state, formed within a short time, and stable in the absence of diffusion for the bacteria. The latter corresponds to a metastable state for small diffusion.

\subsection{The experiment of Swain and Ray}

Finally, we reproduce numerically the experimental observations of Swain and Ray \cite{SwaRay}. To this end we take, as in the experiments, the initial concentration of the fungus to be uniform and constant, with value slightly above the threshold concentration $u_* = 0.2$. Hence, our initial condition was taken as $u(x,y,0) = 0.21$. Now we plant two localized colonies of bacteria in the form 
\begin{equation}
\label{civT3}
v(x,y,0) = 3 \big( e^{-1000((x-0.2)^2 + (y-0.5)^2)} + e^{-1000((x-0.8)^2 + (y-0.5)^2)}\big),
\end{equation}
simulating the administration of the bacteria \textit{B. subtilis} into the control of the \textit{F. oxysporum}. The initial concentration of the chemical was taken as $c(x,y,0) = 0$, as before. In this numerical simulation we set the diffusion coefficient for the bacteria as $D_v = 0$.

\begin{figure}[ht!]
\begin{center}
        \subfigure[$t = 0.1569$.]{
            \label{figT3v-t01569}
            \includegraphics[scale=0.22]{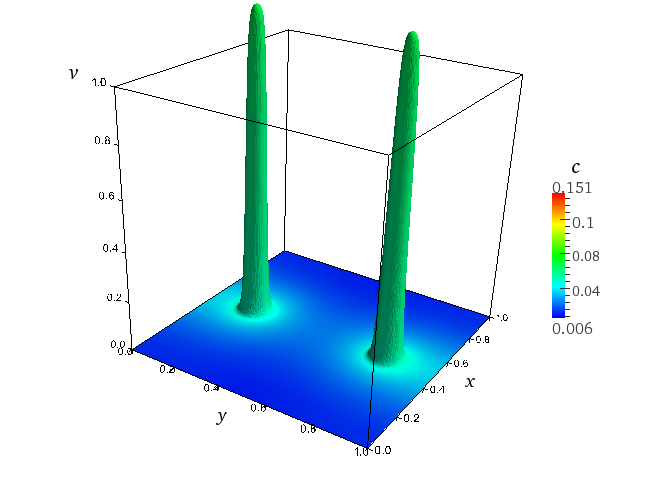}
        }
	\subfigure[$t = 0.5882$.]{
            \label{figT3v-t05882}
            \includegraphics[scale=0.22]{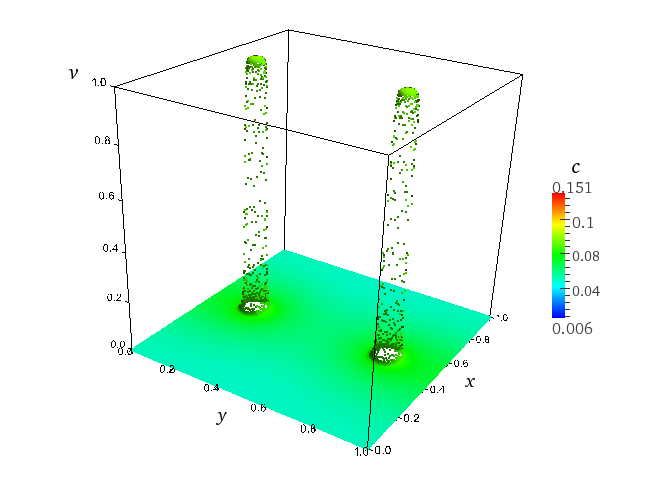}
        }\\
	\subfigure[$t = 0.9804$.]{
            \label{figT3v-t09804}
            \includegraphics[scale=0.22]{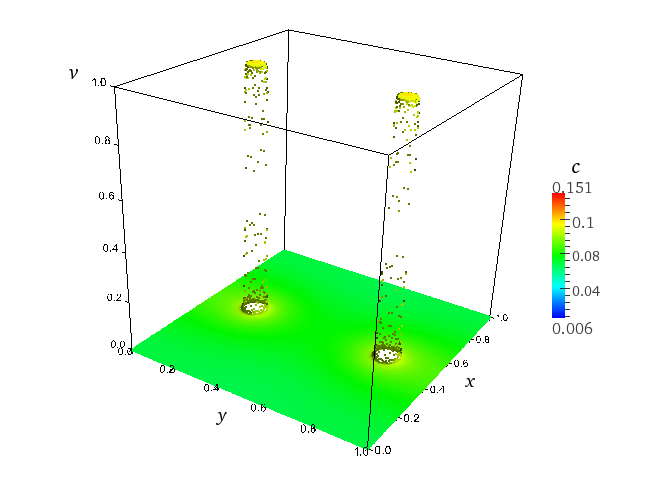}
        }
	\subfigure[$t = 9.8039$.]{
            \label{figT3v-t98039}
            \includegraphics[scale=0.22]{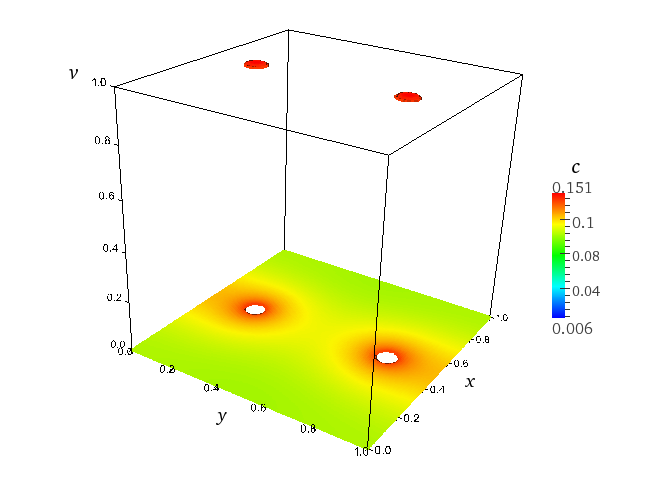}
        }\\
	\subfigure[Chemical concentration $c$ at time $t = 9.8039$.]{
	    \label{figT3c-t98039}
	    \includegraphics[scale=0.22]{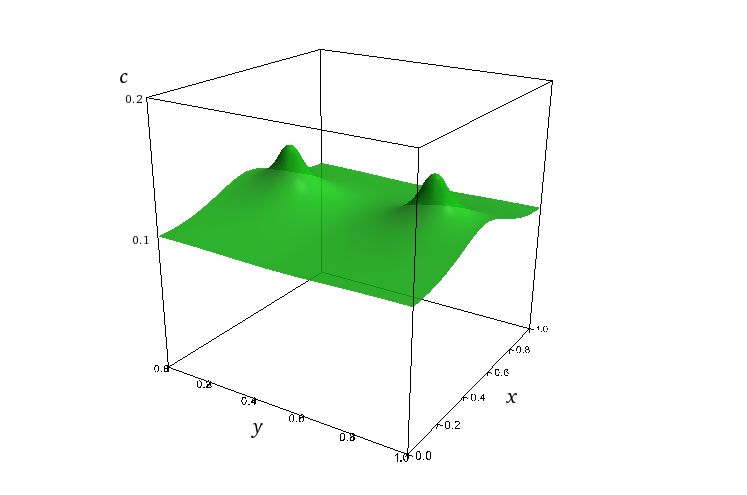}
	}
\end{center}
\caption{\label{figT3v}Concentration $v$ of the bacteria colony for different times. It reaches a stationary state given by the two plateaus in figure (d), and induces, in turn, the steady state for the chemical concentration $c$ shown in figure (e).}
\end{figure}

\begin{figure}[t]
\begin{center}
        \subfigure[$t = 0$.]{
            \label{figT3u-condini}
            \includegraphics[scale=0.19]{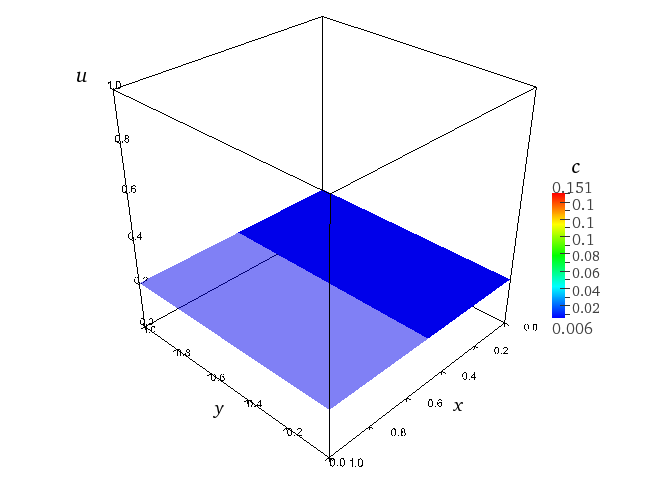}
        }
	\subfigure[$t = 0.0392$.]{
            \includegraphics[scale=0.19]{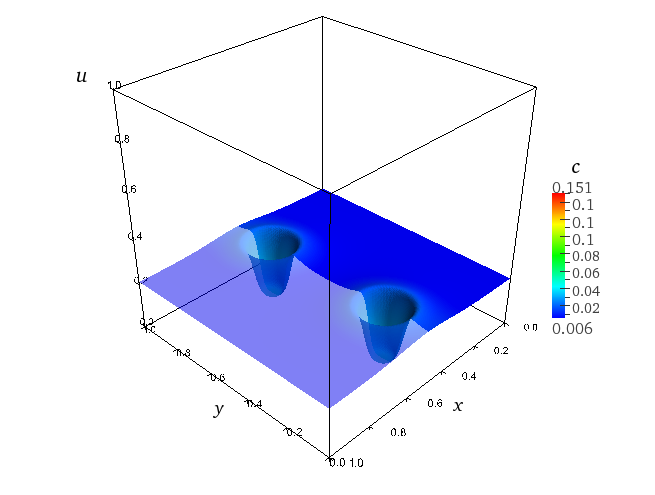}
        }\\
	\subfigure[$t = 0.1176$.]{
            \includegraphics[scale=0.19]{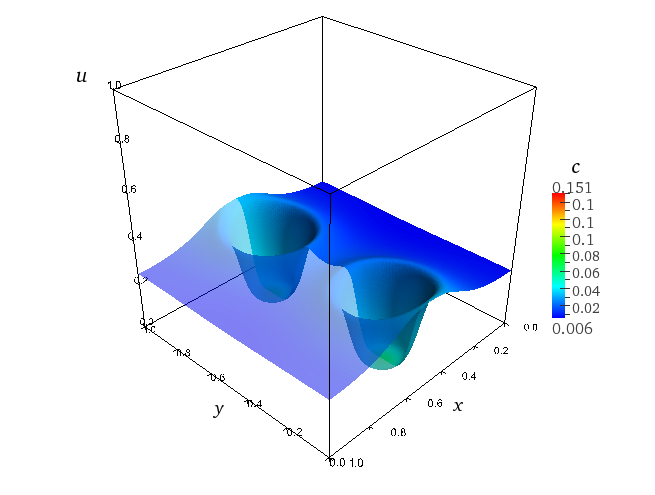}
        }
	\subfigure[$t = 0.1961$.]{
            \includegraphics[scale=0.19]{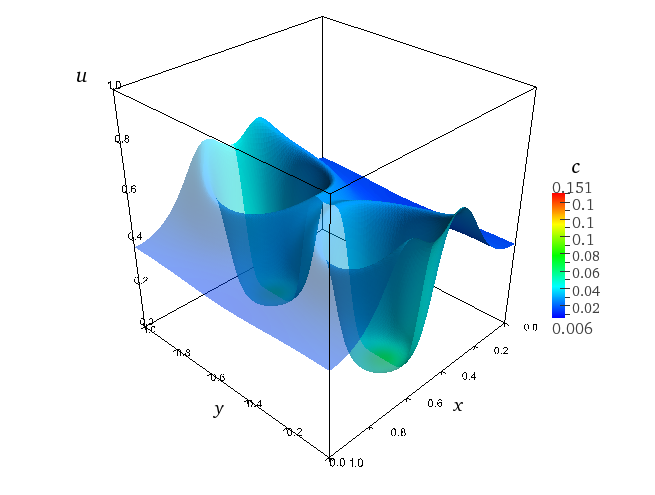}
        }\\
	\subfigure[$t = 0.2549$.]{
            \includegraphics[scale=0.19]{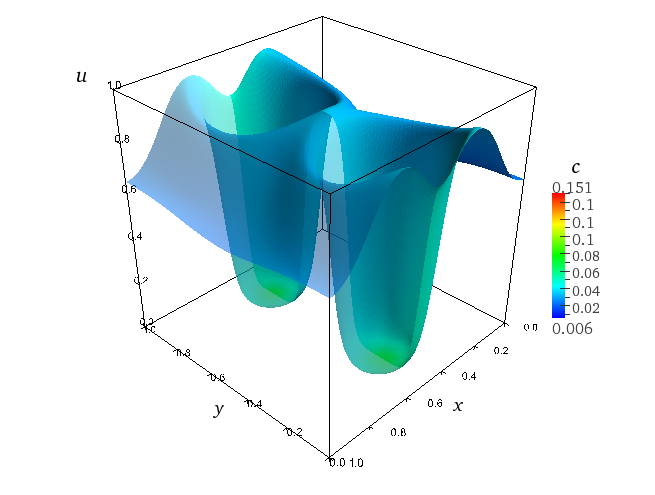}
        }
	\subfigure[$t = 0.3137$.]{
            \includegraphics[scale=0.19]{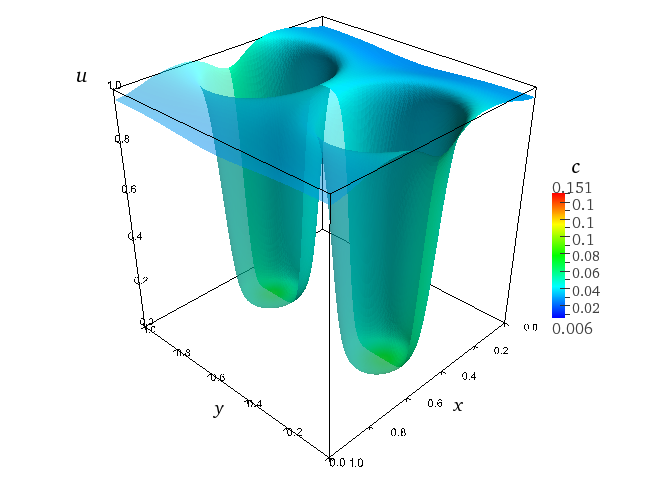}
        }\\
	\subfigure[$t = 0.4902$.]{
            \includegraphics[scale=0.19]{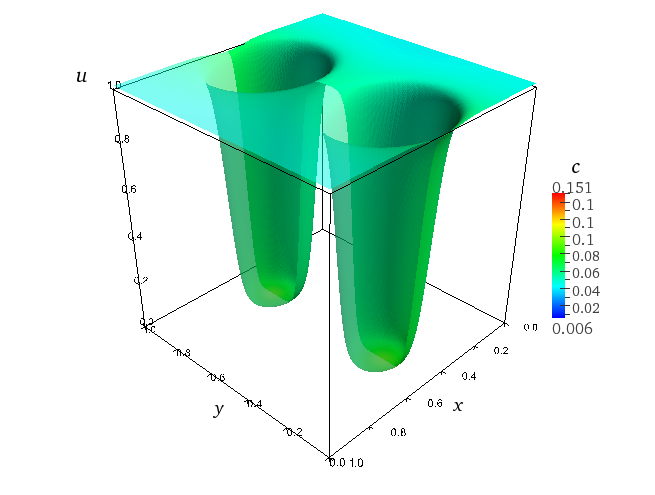}
        }
	\subfigure[$t = 9.8039$.]{
            \label{figT3u-t98039}
            \includegraphics[scale=0.19]{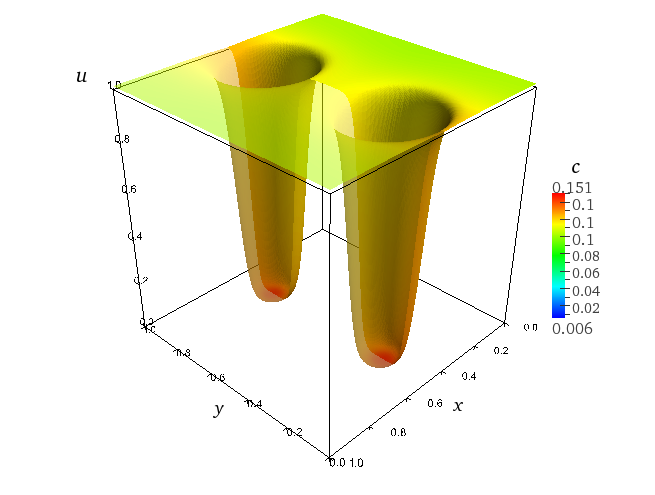}
        }
	\\
	
\end{center}
\caption{\label{figT3u}Concentration $u$ for different times when $D_v = 0$. The initial condition shown in figure (a) is the uniform concentration $u = 0.21$. The fungus diffuses and eventually reaches a stationary state shown in figure (h). Notice that as it evolves, it avoids the repelling region induced by the chemical gradient produced by $v$, whose concentration is shown in figure \ref{figT3c-t98039}. The result is the formation of two equilibrium fronts.}
\end{figure}

Figures \ref{figT3v} and \ref{figT3u} show the dynamics of the solutions under these initial conditions. As expected, the bacteria and the concentration of the chemical evolve practically independently of the two colonies, as shown in figures \ref{figT3v-t01569}, \ref{figT3v-t05882}, \ref{figT3v-t09804}, \ref{figT3v-t98039} and \ref{figT3c-t98039}. Notice the formation of two localized steady plateaus for the bacteria concentration, as shown in figure \ref{figT3v-t98039}. It induces a steady chemical concentration shown in figure \ref{figT3c-t98039}. As in the previous example, the invasive front for the variable $u$ encounters the chemical gradient, localized around two points of the domain. The repulsion causes again the front to go around the higher chemical concentrations depicted in figure \ref{figT3c-t98039}, settling into the superposition (to leading order) of two circular fronts in equilibrium, resembling the one studied in the first example. This evolution is shown in figures \ref{figT3u-condini} through \ref{figT3u-t98039}.

\begin{figure}[ht]
\begin{center}
\includegraphics[scale=0.45]{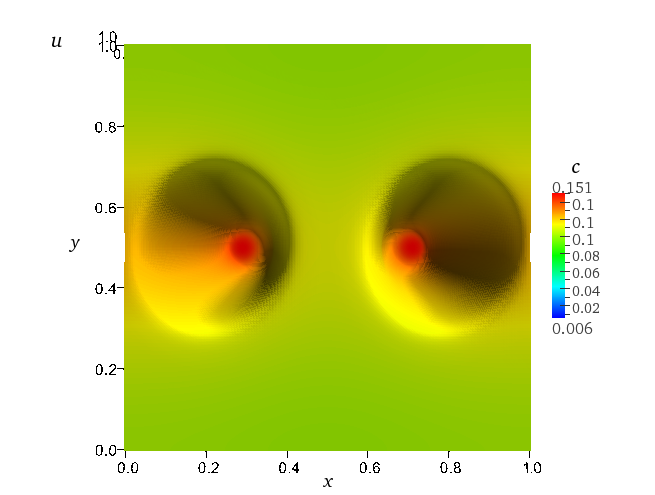}
\end{center}
\caption{\label{figT3utop}Top view of the stationary end state for the concentration $u$ of the harmful colony, at time $t = 9.8039$ (see figure \ref{figT3u-t98039}). It is coloured based on the concentration $c$. Notice the remarkable resemblance with the experimental result of Swain and Ray (figure 2 in \cite{SwaRay}).}
\end{figure}

Finally, figure \ref{figT3utop} shows the top view of the stationary state for the fungus. Notice the remarkable resemblance with the experimental pattern of Swain and Ray (figure 2 in \cite{SwaRay}). The steady state computed here is a good approximation of the metastable pattern which occurs when diffusion of the bacteria is small. We conjecture that the latter is a suitable model for the experimental quasi-stationary state observed in vitro.

\section{Conclusions}
\label{secconclusions}
In this study we have shown how the process of inhibition of an invading front of one species, triggered by the chemical produced by another species, can be accurately described by two reaction-diffusion equations of Nagumo-type which are chemotactically coupled to a third equation for the chemical. We have exhibited the existence of a radially symmetric steady state (for a circular domain), that was shown to be stable in the geometric front propagation limit. It was shown numerically that this steady state is an excellent approximation to the steady state obtained in a square domain as the result of the repulsion by the chemical of an invading front. The asymptotic solutions for the stability of the front suggested instability of the circular steady state as the chemical gradients decrease. This prediction is tested numerically with the expected results. 

In addition, we showed numerically that different colonies of bacteria act independently, producing a final fungus pattern which is well approximated by the superposition of the basic state studied here. This result explains qualitatively the mechanism observed in the in vitro experiments performed by Swain and Ray \cite{SwaRay}. We have thus shown how the solution of a system of equations can be understood in simple terms using the ideas of geometric front propagation, which were able to predict both quantitative and qualitative behaviors. Finally, we anticipate that a distribution of attracting and repelling chemotaxis can produce dendrite-like patterns which can reach a steady state.

\section*{Acknowledgements}

The authors thank Dr. Iv\'an Pavel Moreno for getting them interested in the experiments of Swain and Ray. This research was partially supported by CONACyT (Mexico) and MIUR (Italy), through the MAE Program for Bilateral Research, grant no. 146529.


\bibliography{../refdb/ref}
\bibliographystyle{../refdb/newstyle}


\end{document}